\documentclass[12pt]{amsart}
\usepackage{graphicx}
\usepackage{amssymb}
\usepackage{amsmath}
\usepackage{array}
\usepackage[matrix,frame]{xypic}
\allowdisplaybreaks

\newtheorem{example}{Example}[section]

\newtheorem{theorem}[example]{Theorem}

\newtheorem{definition}[example]{Definition}
\newtheorem{proposition}[example]{Proposition}

\newtheorem{lemma}[example]{Lemma}

\def\Proof{\noindent \it Proof -- \rm}
\def\qed{\hspace{3.5mm} \hfill \vbox{\hrule height 3pt depth 2 pt width 2mm}
\bigskip}

\def\E{{\mathbb E}}
\def\X{{\bf X}}
\def\Y{{\bf Y}}
\def\D{{\bf D}}

\def\Ab{{\bar A}}

\def\B{{\bf B}}

\def\moinsu{{m(\epsilon)}}

\def\FQSym{{\bf FQSym}}

\def\Sym{{\bf Sym}}
\def\NCSF{{\bf Sym}}

\def\PBT{{\bf PBT}}
\def\CBT{{\bf CBT}}

\def\Std{{\rm Std}}
\def\cstd{{\bf Std}}

\def\convol{{*}}

\def\<{\langle}
\def\>{\rangle}

\def\C{\operatorname{\mathbb C}}

\def\Z{\operatorname{\mathbb Z}}

\def\G{{\bf G}}

\def\SG{{\mathfrak S}}

\def\A{{\bf A}}

\def\Des{\operatorname{Des}}
\def\Bdes{\operatorname{Bdes}}

\newcommand{\free}[1]{\langle#1\rangle}

\def\MR{{\bf MR}}

\def\BSym{{\bf BSym}}

\def\std{{\rm std}}
\def\ol{\overline}

\usepackage{boxedminipage}
\usepackage[usenames]{color}
\newenvironment{NOTE}[1][NOTE]
 {\bigskip\begin{center}\color{red}\begin{boxedminipage}{4.5in}\setlength{\parindent}{1em}\noindent\color{blue}\textbf{#1. }}
 {\end{boxedminipage}\end{center}\bigskip}

\def\gf#1#2{\genfrac{}{}{0pt}{}{#1}{#2}}

\def\SEC{{\bf sec}}
\def\TAN{{\bf tan}}
\def\COS{{\bf cos}}
\def\SIN{{\bf sin}}

\def\id{{\rm id}}
\def\sSym{{\bf sSym}}
\def\Sg{{\sf Sg}}
\def\II{{C}}

\title{The algebraic combinatorics of snakes}

\author[M. Josuat-Verg\`es, J.-C. Novelli, and J.-Y. Thibon]%
{Matthieu Josuat-Verg\`es, Jean-Christophe Novelli,\\
and Jean-Yves Thibon}

\thanks{M. Josuat-Verg\`es was supported by the Austrian Science Foundation
(FWF) via the grant Y463.}
\address{
Fak\"ult\"at f\"ur Mathematik \\
Universit\"at Wien \\
Garnisongasse 3, 1090 Wien \\
Austria}
\address{
Institut Gaspard Monge, Universit\'e Paris-Est Marne-la-Vall\'ee \\
5 Boulevard Des\-cartes \\
Champs-sur-Marne \\
77454 Marne-la-Vall\'ee cedex 2 \\ France}

\email[Matthieu Josuat-Verg\`es]{Matthieu.Josuat-Verges@univie.ac.at}
\email[Jean-Christophe Novelli]{novelli@univ-mlv.fr}
\email[Jean-Yves Thibon]{jyt@univ-mlv.fr} 
\date{\today}

\begin{document}

\begin{abstract}
Snakes are analogues of alternating permutations defined for any Coxeter
group. We study these objects from the point of view of combinatorial Hopf
algebras, such as noncommutative symmetric functions and their
generalizations.
The main purpose is to show that several properties of the generating
functions of snakes, such as differential equations or closed form as
trigonometric functions, can be lifted at the level of noncommutative
symmetric functions or free quasi-symmetric functions. The results take the
form of algebraic identities for type B noncommutative symmetric functions,
noncommutative supersymmetric functions and colored free quasi-symmetric
functions.
\end{abstract}

\maketitle

\section{Introduction}

Snakes, a term coined by Arnol'd \cite{Arnold}, are generalizations of
alternating permutations.
These permutations arose as the solution of what is perhaps the first
example of an inverse problem in the theory of generating functions: given
a function whose Taylor series has nonnegative integer coefficients, find a
family of combinatorial objects counted by those coefficients. For example,
in the expansions
\begin{equation}
\tan z = \sum_{n\ge 0} E_{2n+1}\frac{z^{2n+1}}{(2n+1)!}\quad
\text{and}
\quad
\sec z =  \sum_{n\ge 0} E_{2n}\frac{z^{2n}}{(2n)!}\,,
\end{equation}
the coefficients $E_n$ are nonnegative integers. 

It was found in 1881 by D. Andr\'e \cite{Andre} that $E_n$ was the number of
alternating permutations in the symmetric group $\SG_n$.

Whilst this result is not particularly difficult and can be proved in several
ways, the following explanation is probably not far from being optimal: there
exists an associative (and noncommutative) algebra admitting a basis labelled
by all permutations, and such that the map $\phi$ sending any $\sigma\in\SG_n$
to $\frac{z^n}{n!}$ is a homomorphism.
In this algebra, the formal series 
\begin{equation}
C = \sum_{n\ge 0}(-1)^n \id_{2n}\quad \text{and}\quad
S = \sum_{n\ge 0}(-1)^n \id_{2n+1}
\end{equation}
(alternating sums of even and odd identity permutations) are respectively
mapped to $\cos z$ and $\sin z$ by $\phi$. The series $C$ is clearly
invertible, and one can see by a direct calculation that
$C^{-1} + C\cdot S^{-1}$ is the sum of all alternating
permutations~\cite{NCSF1}.

Such a proof is not only illuminating, it says much more than the original
statement. For example, one can now replace $\phi$ by more complicated
morphisms, and obtain generating functions for various statistics on
alternating permutations.

The symmetric group is a Coxeter group, and snakes are generalizations of
alternating permutations to arbitrary Coxeter groups. Such generalizations
were first introduced by Springer \cite{Springer}. For the infinite series
$A_n$, $B_n$, $D_n$, Arnol'd \cite{Arnold} related the snakes to the geometry
of bifurcation diagrams.

The aim of this article is to study the snakes of the classical Weyl groups
(types A, B and D) by noncommutative methods, and to generalize the results
to some series of wreath products (colored permutations).

The case of symmetric groups (type A) is settled by the algebra of Free
quasi-symmetric functions $\FQSym$ (also known the Malvenuto-Reutenauer
algebra) which is based on permutations,  and its subalgebra $\Sym$
(noncommutative symmetric functions), based on integer compositions. To deal
with the other types, we need an algebra based on signed permutations, and
some of its subalgebras defined by means of the superization map introduced
in~\cite{NT-super}.

After reviewing the necessary background and the above mentioned proof of
the result of Andr\'e, we recover results of Chow~\cite{Chow} on type B snakes,
and derive some new generating functions for this type. This suggests a
variant of the definition of snakes, for which the noncommutative generating
series is simpler.
These considerations lead us to some new identities satisfied by the
superization map on noncommutative symmetric functions.  Finally, we propose
a completely different combinatorial model for the generating function of type
B snakes, based on  interesting identities in the algebra of signed
permutations. We also present generalizations of (Arnol'd's) Euler-Bernoulli
triangle, counting alternating permutations according to their last value, and
extend the results to wreath products and to type D, for which we propose an
alternative definition of snakes.

\section{Permutations and noncommutative trigonometry}

\subsection{Free quasi-symmetric functions}

The simplest way to define our algebra based on permutations is by means of
the classical \emph{standardization process}, familiar in combinatorics and in
computer science.
Let $A=\{a_1,a_2,\dots\}$ be an infinite totally ordered alphabet.
The \emph{standardized word} $\Std(w)$ of a word $w\in A^*$ is the permutation
obtained by iteratively scanning $w$ from left to right, and labelling
$1,2,\dots$ the occurrences of its smallest letter, then numbering the
occurrences of the next one, and so on.
Alternatively, $\sigma=\std(w)^{-1}$ can be characterized as the unique
permutation of minimal length such that $w\sigma$ is a nondecreasing word.
For example, $\std(bbacab)=341625$.

We can now define polynomials 
  \begin{equation}
    \G_{\sigma}(A) :=
    \sum_{\std(w) = \sigma} w\,.
  \end{equation}
It is not hard to check that these polynomials span a subalgebra of $\C\<A\>$,
denoted by $\FQSym(A)$, an acronym for {\em Free Quasi-Symmetric functions}.

The multiplication rule is, for $\alpha\in\SG_k$ and
$\beta\in\SG_\ell$,
\begin{equation}
\label{multG}
\G_\alpha \G_\beta
=
\sum_{\gamma\in\alpha\convol\beta} \G_\gamma,
\end{equation}
where $\alpha*\beta$ is the set of permutations $\gamma\in\SG_{k+\ell}$ such
that $\gamma=u\cdot v$ with $\std(u)=\alpha$ and $\std(v)=\beta$. This is the
\emph{convolution} of permutations (see~\cite{Reu}).
%
Note that the number of terms in this product depends only on $k$ and $\ell$,
and is equal to the binomial coefficient $\binom{k+\ell}{k}$. Hence, the
map
\begin{equation}
\phi:\ \sigma\in\SG_n \longmapsto \frac{z^n}{n!}
\end{equation}
is a homomorphism of algebras $\FQSym\rightarrow \C[z]$.

\subsection{Noncommutative symmetric functions}

The algebra $\Sym(A)$ of noncommutative symmetric functions  over $A$
is the subalgebra of $\FQSym$ generated by the
identity permutations \cite{NCSF1,NCSF6}
\begin{equation}
S_n(A) := \G_{12\dots n}(A)
 = \sum_{i_1\le i_2\le\dots\le i_n}a_{i_1}a_{i_2}\dots a_{i_n}\,.
\end{equation}
These polynomials are obviously algebraically independent, so that the products
\begin{equation}
S^I := S_{i_1}S_{i_2}\dots S_{i_r}
\end{equation}
where $I=(i_1,i_2,\dots,i_r)$ runs over compositions of $n$, form a basis of
$\Sym_n$, the homogeneous component of degree $n$ of $\Sym$.

Recall that a descent of a word $w=w_1w_2\dots w_n\in A^n$ is an index $i$
such that $w_i>w_{i+1}$. The set of such $i$ is denoted by $\Des(w)$.
Hence, $S_n(A)$ is the sum of all nondecreasing words of length $n$ (no
descent), and $S^I(A)$ is the sum of all words which may have a descents only
at places from the set
\begin{equation}
\Des(I) = \{ i_1,\ i_1+i_2, \dots , i_1+\dots+i_{r-1}\},
\end{equation}
called the descent set of $I$. Another important basis is
\begin{equation}
\label{ncribbon}
R_I(A) = \sum_{\Des(w)=\Des(I)}w = \sum_{\Des(\sigma)=\Des(I)}\G_\sigma\,,
\end{equation}
the ribbon basis, formed by sums of words having descents exactly at
prescribed places. From this definition, it is obvious that
if $I=(i_1,\dots,i_r)$, $J=(j_1,\dots,j_s)$
\begin{equation}
R_I(A) R_J(A) = R_{IJ}(A) + R_{I\triangleright J}(A)
\end{equation}
with $IJ=(i_1,\dots,i_r,j_1,\dots,j_s)$ and
$I\triangleright J=(i_1,\dots,i_r+j_1,j_2,\dots,j_s)$.

\subsection{Operations on alphabets}

If $B$ is another totally ordered alphabet, we denote by $A+B$ the ordinal sum
of $A$ and $B$. This allows to define noncommutative symmetric functions of
$A+B$, and
\begin{equation}
S_n(A+B)=\sum_{i=0}^nS_i(A)S_{n-i}(B) \quad (S_0=1).
\end{equation}
If we assume that $A$ and $B$ commute, this operation defines a coproduct, for
which $\Sym$ is a graded bialgebra, hence a Hopf algebra. The same is true of
$\FQSym$.
Symmetric functions of the virtual alphabet $(-A)$ are defined by the condition
\begin{equation}
\sum_{n\ge 0}S_n(-A)=\left(\sum_{n\ge 0}S_n(A)\right)^{-1}
\end{equation}
and more generally, for a difference $A-B$,
\begin{equation}
\sum_{n\ge 0}S_n(A-B)=\left(\sum_{k\ge 0}S_k(B)\right)^{-1}\sum_{l\ge 0}S_l(A)
\end{equation}
(note the reversed order, see \cite{NCSF2} for detailed explanations).

\subsection{Noncommutative trigonometry}

\subsubsection{Andr\'e's theorem}

One can now define ``noncommutative trigonometric functions'' by
\begin{equation}
\COS(A)=\sum_{n\ge 0}(-1)^nS_{2n}(A)\ \text{\quad and\quad}\
\SIN(A)= \sum_{n\ge 0}(-1)^nS_{2n+1}(A).
\end{equation}
The image by $\phi$ of these series are the usual trigonometric functions.
With the help of the product formula for the ribbon basis,
it is easy to see that
\begin{equation}
\SEC := \COS^{-1} =\sum_{n\ge 0}R_{(2^n)}\ \text{\quad and\quad }\
\TAN := \COS^{-1}\SIN  = \sum_{n\ge 0}R_{(2^n1)}
\end{equation}
which implies Andr\'e's theorem: the coefficient of $\frac{z^n}{n!}$ in
$\sec(z)+\tan(z)$ is the number of alternating permutations of $\SG_n$
(if we choose to define alternating permutations as those of shape
$(2^n)$ and $(2^n1)$). In $\FQSym$,
\begin{equation}
\SEC + \TAN= \sum_{\sigma\ \text{alternating}}\G_\sigma\,.
\end{equation}

\subsubsection{Differential equations}

If $\partial$ is the derivation of $\Sym$ such that $\partial S_n=S_{n-1}$,
then
\begin{equation}
X=\TAN =\sum_{m\ge 0}R_{(2^m1)}
\ \text{ and }\ Y=\SEC=\sum_{m\ge 0} R_{(2^m)}
\end{equation}
satisfy the differential equations
\begin{equation}
\partial X=1+X^2,\quad \partial Y=XY\,.
\end{equation}
These equations can be lifted to $\FQSym$, actually to its
subalgebra $\PBT$, the Loday-Ronco algebra of planar binary trees
(see~\cite{HNT} for details). Solving them in this algebra
provides yet another combinatorial proof of Andr\'e's result.

Let us sketch it for the tangent.
The original proof of Andr\'e relied upon the differential equation
\begin{equation}
\frac{dx}{dt} = 1 + x^2
\end{equation}
whose $x(t)=\tan(t)$ is the solution such that $x(0)=0$. Equivalently,
$x(t)$ is the unique solution of the functional equation
\begin{equation}
x(t)= t + \int_0^t x(s)^2 ds
\end{equation}
which can be solved by iterated substitution.

In general, given an associative algebra $R$, we can consider the functional
equation for the power series $x\in R[[t]]$
\begin{equation}
\label{eqbin}
x = a + B(x,x)
\end{equation}
where $a\in R$ and $B(x,y)$ is a bilinear map with values in $R[[t]]$, such
that the valuation of $B(x,y)$ is strictly greater than the sum of the
valuations of $x$ and $y$. Then, Equation~(\ref{eqbin}) has a unique solution
\begin{equation}
\label{solbin}
x = a + B(a,a) + B(B(a,a),a)+ B(a,B(a,a))+ \dots =  \sum_{T\in\CBT} B_T(a)
\end{equation}
where $\CBT$ is the set of (complete) binary trees, and for a tree $T$,
$B_T(a)$ is the result of evaluating the expression formed by labeling
by $a$ the leaves of $T$ and by $B$ its internal nodes.
Pictorially,
  \begin{gather*}
    x = a + B(a,a) + B(B(a,a), a) + B(a, B(a,a)) + \dots \\
    = a +
    \vcenter{\xymatrix@C=2mm@R=2mm{
        *{} & {B}\ar@{-}[dr]\ar@{-}[dl] \\
        a   & *{} & a
      }} +
    \vcenter{\xymatrix@C=2mm@R=2mm{
        *{} & *{} & {B}\ar@{-}[dr]\ar@{-}[dl] \\
        *{}   & {B}\ar@{-}[dr]\ar@{-}[dl] & *{} & a \\
        a   & *{} & a
      }} +
    \vcenter{\xymatrix@C=2mm@R=2mm{
        *{} & {B}\ar@{-}[dr]\ar@{-}[dl] \\
        a & *{} & {B}\ar@{-}[dr]\ar@{-}[dl] \\
        *{} & a   & *{} & a
      }} + \dots
  \end{gather*}

It is proved in \cite{HNT} that if one defines
\begin{equation}
\label{part}
\partial \G_\sigma := \G_{\sigma'}\,
\end{equation}
where $\sigma'$ is obtained from $\sigma$ by erasing its maximal letter $n$,
then $\partial$ is a derivation of $\FQSym$. Its restriction to $\Sym$
coincides obviously with the previous definition.
For $\alpha\in\SG_k$, $\beta\in\SG_\ell$, and $n=k+\ell$, set
\begin{equation}
\B(\G_\alpha,\G_\beta)
=\sum_{\gf{\gamma=u(n+1)v}{\std(u)=\alpha,\std(v)=\beta}}\G_\gamma\,.
\end{equation}
Clearly,
\begin{equation}
\partial \B(\G_\alpha,\G_\beta)=\G_\alpha\G_\beta\,,
\end{equation}
and our differential equation for the noncommutative tangent is now replaced
by the fixed point problem
\begin{equation}
\label{eqbinX}
X = \G_1 + \B(X,X)\,.
\end{equation}
of which it is the unique solution. Again, solving it by iterations gives back
the sum of alternating permutations. As an element of the Loday-Ronco algebra,
$\TAN$ appears as the sum of all permutations whose decreasing tree is complete.

The same kind of equation holds for $Y$:
\begin{equation}
\label{eqbinY}
Y = 1 + \B(X,Y)\,.
\end{equation}

Hence, the noncommutative secant $\SEC$ is therefore an element of $\PBT$, so
a sum of binary trees.
The trees are well-known: they correspond to complete binary trees (of odd
size) where one has removed the last leaf.

\subsection{Derivative polynomials}

For the ordinary tangent and secant, the differential equations
imply the existence~\cite{Hof} of two sequences of polynomials $P_n$, $Q_n$
such that
\begin{equation}
\frac{d^n}{dz^n}(\tan z) = P_n(\tan z)
\quad\text{and}\quad
\frac{d^n}{dz^n}(\sec z) = Q_n(\tan z)\sec z.
\end{equation}
Since $\partial$ is a derivation of $\Sym$, we have as well for the
noncommutative lifts
\begin{equation}
\partial^n(X)=P_n(X)\quad\text{and}\quad \partial^n(Y)=Q_n(X)Y\,.
\end{equation}
Hoffman~\cite{Hof} gives the exponential generating functions
\begin{equation}
P(u,t)=\sum_{n\ge 0}P_n(u)\frac{t^n}{n!}
=\frac{\sin  t + u\cos t}{\cos t - u\sin t}\ \text{and}\
Q(u,t)=\sum_{n\ge 0}Q_n(u)\frac{t^n}{n!}
=\frac{1}{\cos t - u\sin t}\,.
\end{equation}
The noncommutative version of these identities can be readily derived as
follows. We want to compute
\begin{equation}
P(X,t)=\sum_{n\ge 0}P_n(X)\frac{t^n}{n!}=e^{t\partial}X\,.
\end{equation}
Since $\partial$ is a derivation, $e^{t\partial}
$ is an automorphism of
$\Sym$. It acts on the generators $S_n$ by
\begin{equation}
e^{t\partial}S_n(A)=\sum_{k=0}^nS_{n-k}(A)\frac{t^k}{k!} = S_n(A+t\E)
\end{equation}
where $t\E$ is the ``virtual alphabet'' such that $S_n(t\E)=\frac{t^n}{n!}$.
Hence,
\begin{equation}
\begin{split}
P(X,t)& =\TAN(A+t\E)=\COS(A+t\E)^{-1}\SIN(A+t\E) \\
&=(\cos t -X\sin t)^{-1}(\sin t +X\cos t)
\end{split}
\end{equation}
as expected. Similarly,
\begin{align}
  Q_n(X,t) &= \sum_{n\geq 0} Q_n(X)\frac{t^n}{n!} = (e^{t\partial} Y) Y^{-1} \\
               &= \COS(A+t\E)^{-1} \COS(A)= (\cos t - X \sin t)^{-1}.
\end{align}

\section{The uniform definition of snakes for Coxeter groups}
\label{section2}

Before introducing the relevant generalizations of $\Sym$ and
$\FQSym$, we shall comment on the definitions of snakes and alternating
permutations for general Coxeter groups.

It is apparent in Springer's article \cite{Springer} that alternating 
permutations can be defined in a uniform way for any Coxeter group.
Still, little attention has been given to this fact. For example, 20 years
later, Arnol'd~\cite{Arnold} gives separately the definitions of snakes of
type $A$, $B$ and $D$, even though there is no doubt that he was aware of the
uniform definition.
The goal of this section is to give some precisions and to simplify the proof of
Springer's result in \cite{Springer}.

Let $(W,S)$ be an irreducible Coxeter system. Recall that $s\in S$ is a
descent of $w\in W$ if $\ell(ws)<\ell(w)$ where $\ell$ is the length function.
When $J\subset S$, we denote by $\mathcal{D}_J$ the \emph{descent class}
defined as
$\{ w\in W \; : \; \ell(ws)<\ell(w) \Leftrightarrow s\in J \} $.
Following Arnol'd, let us consider the following definition.
\begin{definition}
The Springer number of the Coxeter system $(W,S)$ is
\begin{equation}
\label{def_K_des}
  K(W) := \max_{J\subseteq S} (\# \mathcal{D}_J).
\end{equation}
\end{definition}
The aim of \cite{Springer} is to give a precise description of sets $J$
realizing this maximum. The result is as follows:

\begin{theorem}[Springer \cite{Springer}]
\label{subsets}
Let $J\subseteq S$. Then, $ K(W) =  \# \mathcal{D}_J $ if and only if $J$
and $S\backslash J$ are independent subsets of the Coxeter graph $S$
(\emph{i.e.}, they contain no two adjacent vertices).
In particular, there are two such subsets $J$ which are complementary to each
other.
\end{theorem}

Therefore we can choose a subset $J$ such that $ K(W) = \#  \mathcal{D}_J $
and call {\it snakes} of $W$ the elements of $\mathcal{D}_J$. The other choice
$S\backslash J$ would essentially lead  to the same objects, since there is a
simple involution on $W$ exchanging the subsets $\mathcal{D}_J$ and
$\mathcal{D}_{S\backslash J}$.

We are thus led to the following definition:

\begin{definition}
Let $(W,S)$ be a Coxeter group, and $J$ be a maximal independent subset of $S$.
The \emph{snakes} of $(W,S)$ are the elements of the descent class
$\mathcal{D}_J$.
\end{definition}

This definition depends on the choice of $J$, so that we can consider two
families of snakes for each $W$. In the case of alternating permutations,
these are usually called the up-down and down-up permutations, and are
respectively defined by the conditions
\begin{equation}
\sigma_1<\sigma_2>\sigma_3<\dots \text{\quad or\quad}
\sigma_1>\sigma_2<\sigma_3>\dots
\end{equation}

It is natural to endow a descent class with the restriction of the weak order,
and this defines what we can call the \emph{snake poset} of $(W,S)$. Known
results show that this poset is a lattice~\cite{Bj}.
Let us now say a few words about the proof of Theorem~\ref{subsets}, which
relies upon the following lemma.

\begin{lemma}[Springer \cite{Springer}]
\label{lemspringer}
Let $J\subseteq S$, and assume that there is an edge $e$ of the Coxeter graph
whose endpoints are both in $J$ or both not in $J$. Let $S=S_1\cup S_2$ be the
connected components obtained after removing $e$.
Let $J' = (S_1\cap J) \cup (  S_2 \cap (S\backslash J))$.
Then, $\#\mathcal{D}_{J} < \#\mathcal{D}_{J'}$.
\end{lemma}

Using the above lemma, we see that if $J$ or complementary is not independent,
we can find another subset $J'$ having a strictly bigger descent class, and
Theorem~\ref{subsets} then follows.
Whereas Lemma~\ref{lemspringer} is the last one out of a series of 5 lemmas in
Springer's article, we give here a simple geometric argument.

Let $R$ be a root system for $(W,S)$, and let $\Pi = (\alpha_s)_{s\in S}$ be a
set of simple roots. There is a bijection $w\mapsto i(w)$ between $W$ and the
set of Weyl chambers, so that $s\in S$ is a descent of $w\in W$ if and only if
$i(w)$ lies in the half-space
$\{ v\in\mathbb{R}^n : \langle v,\alpha_s\rangle \geq 0 \}$.
For any $J\subset S$, let
\begin{equation}
    C_J = \{  v\in\mathbb{R}^n :  \langle v,\alpha_s\rangle \geq 0
          \hbox{ if } s\in J, \hbox{ and }
          \langle v,\alpha_s\rangle \leq 0 \hbox{ if } s\notin J    \}.
\end{equation}
It is the closure of the union of Weyl chambers $i(w)$ where
$w\in\mathcal{D}_J$.
Now, let $e$, $S_1$, $S_2$ and $J'$ be as in the lemma, and let $x\in S_1$ and
$y\in S_2$ be the endpoints of $e$.  Note that either $x$ or $y$ is in $J'$
but not both.
Let $\sigma$ be the orthogonal symmetry through the linear span of
$\{ \alpha_j : j\in S_1 \}$.
We claim that $\sigma( C_J )\subsetneq C_{J'}$, and this implies
$\#\mathcal{D}_{J} < \#\mathcal{D}_{J'}$ since $C_{J}$
contains strictly less Weyl chambers than $C_{J'}$.

So it remains to show that $\sigma( C_J )\subsetneq C_{J'}$. It is convenient
to use the notion of {\it dual cone}, which is defined for any closed convex
cone $C\subset \mathbb{R}^n$ as
\begin{equation}
  C^* := \{ v\in\mathbb{R}^n \; : \;  \langle v , w \rangle \geq 0
          \hbox{ for any } w\in C \}.
\end{equation}
The map $C\mapsto C^*$ is an inclusion-reversing involution on closed convex
cones, and it commutes with any linear isometry, so that we have to prove that
$\sigma( C_{J'}^* ) \subsetneq C_{J}^*$.
Since $(C_1\cap C_2)^*= C_1^*+C_2^*$ and since the dual of the
half-space $\{v\in\mathbb{R}^n\; : \; \langle v,w \rangle\geq 0 \}$ is the
half-line $\mathbb{R}^+w$, the dual of $C_{J'}$ is
\begin{equation}
   C_{J'}^* =  \bigg\{  \sum_{s\in S} u_s \alpha_s \; : \;  u_s \geq 0
   \hbox{ if } s\in J', \;  \hbox{ and } u_s \leq 0 \hbox{ if } s\notin J'
               \bigg\},
\end{equation}
and the same holds for $C_J^*$.
A description of $\sigma(C_{J'}^*)$ is obtained easily by linearity since
$\sigma(\alpha_s)=\alpha_s$ if $s\in S_1$,
$\sigma(\alpha_s)=-\alpha_s$ if $s\in S_2\backslash y$, and
$\sigma(\alpha_y)
= -\alpha_y + 2 \langle \alpha_x,\alpha_y \rangle
   \langle \alpha_x,\alpha_x \rangle ^{-1} \alpha_x$.
Indeed, let $v = \sum_{s\in S} u_s \alpha_s \in C_{J'}^*$, we have:
\begin{equation}
   \sigma(v)  =
   \sum_{s\in S_1} u_s \alpha_s
 - \sum_{s\in S_2\backslash y} u_s \alpha_s
 - u_y \alpha_y + 2 u_y \langle \alpha_x,\alpha_y \rangle
                        \langle \alpha_x,\alpha_x \rangle ^{-1} \alpha_x.
\end{equation}
Since $u_x$ and $u_y$ have different signs and
$\langle \alpha_x,\alpha_y \rangle<0$, we obtain $\sigma(v)\in C_J^*$.
We have thus proved $\sigma( C_{J'}^* ) \subset C_{J}^*$. To show the strict
inclusion, note that either $\alpha_y$ or $-\alpha_y$ is in $C_J^*$. But none
of these two elements is in $\sigma( C_{J'}^* )$, because in the above
formula for $\sigma(v)$, if $u_y\neq 0$ there is a nonzero term in $\alpha_x$
as well. This completes the proof.

\section{Signed permutations and combinatorial Hopf algebras}

Whereas the constructions of the Hopf algebras $\NCSF$, $\PBT$, and $\FQSym$
appearing when computing the usual tangent are almost straightforward, the
situation is quite different in type $B$. First, there are at least three
different generalizations of $\NCSF$ to a pair of alphabets, each with its
own qualities either combinatorial or algebraic.
Moreover, there are also two different ways to generalize $\FQSym$. The
generalizations of $\PBT$ are not (yet) defined in the literature but the
computations done in the present paper give a glimpse of what they should be.

Here follows how they embed each in the other. All embeddings are embeddings
of Hopf algebras except the two embeddings concerning $\BSym$ which is not
itself an algebra. However, the embedding of $\NCSF(A|\Ab)$ into $\NCSF^{(2)}$
obtained by composing the two previous embeddings is a Hopf embedding:

\begin{equation}
\NCSF(A|\Ab) \hookrightarrow \BSym \hookrightarrow \NCSF^{(2)}
\hookrightarrow \FQSym(A|\Ab) \hookrightarrow \FQSym^{(2)}.
\end{equation}


\subsection{The Mantaci-Reutenauer algebra of type $B$}

The most straightforward definition of $\NCSF$ in type $B$ is to generalize
the combinatorial objects involved in the definition: change compositions into
signed compositions.

We denote by $\Sym^{(2)}=\MR$~\cite{MR} the free product $\Sym\star\Sym$ of
two copies of the Hopf algebra of noncommutative symmetric functions. In other
words, $\MR$ is the free associative algebra on two sequences $(S_n)$ and
$(S_{\bar n})$ ($n\ge 1$).
We regard the two copies of $\Sym$ as noncommutative symmetric functions on
two auxiliary alphabets: $S_n=S_n(A)$ and $S_{\bar n}=S_n(\Ab)$.
We denote by $F\mapsto \bar F$ the involutive automorphism%
\footnote{This differs from the convention used in some references.}
which exchanges $S_n$ and $S_{\bar n}$.
And we denote the generators of $\NCSF^{(2)}$ by $S_{(k,\epsilon)}$ where
$\epsilon=\{\pm1\}$, so that $S_{(k,1)}=S_k$ and $S_{(k,-1)}= S_{\bar k}$.

\subsection{Noncommutative supersymmetric functions}

The second generalization of $\NCSF$ comes from the transformation of
alphabets sending $A$ to a combination of $A$ and $\Ab$. It is the algebra
containing the type $B$ alternating permutations.

We define $\NCSF(A|\Ab)$ as the subalgebra of $\NCSF^{(2)}$ generated by
the $S_n^\#$ where, for any $F\in \Sym(A)$,
\begin{equation}
\label{diese}
F^\# = F(A|\bar A) = F(A-q\bar A)|_{q=-1}\,,
\end{equation}
called the supersymmetric version, or superization, of $F$ \cite{NT-super}.

The expansion of an element of $\NCSF(A|\Ab)$ as a linear combination in
$\NCSF^{(2)}$ is done thanks to generating series. Indeed,
\begin{equation}
\sigma_1^\# = 
\bar\lambda_1\sigma_1
= \left( \sum_{k\geq0} \overline{\Lambda_k}\right)
  \left( \sum_{m\geq0} S_m \right)
\end{equation}
where $\overline{\Lambda_k} = \sum_{I\models k} (-1)^{\ell(I)-k} \overline{S^I}$,
as follows from $\bar \lambda_1 = ( \bar \sigma_{-1} )^{-1}$
(see~\cite{NCSF2}). For example,

\begin{equation}
S_1^\# = S^1 + S^{\ol1},
\qquad
S_2^\# = S^2 + S^{{\ol1}1} - S^{\ol2} + S^{\ol1\ol1},
\end{equation}
\begin{equation}
S_3^\# = S^3 + S^{\ol12} + S^{\ol1\ol11} - S^{\ol21} + S^{\ol1\ol1\ol1}
         - S^{\ol2\ol1}  - S^{\ol1\ol2} + S^{\ol3}.
\end{equation}

\subsection{Noncommutative symmetric functions of type $B$}

The third generalization of $\NCSF$ is not an algebra but only a cogebra but
is the generalization one gets with respect to the group $B_n$: its graded
dimension is $2^n$ and, as we shall see later in this paragraph, a basis of
$\BSym$ is given by sums of permutations having given descents
\emph{in the type $B$ sense}. This algebra contains the snakes of type $B$.

Noncommutative symmetric functions of type $B$ were introduced in~\cite{Chow}
as the right $\Sym$-module $\BSym$ freely generated by another sequence
$(\tilde S_n)$ ($n\ge 0$, $\tilde S_0=1$) of homogeneous elements, with
$\tilde\sigma_1$ grouplike. This is a coalgebra, but not an algebra.

We embed $\BSym$ as a sub-coalgebra and right sub-$\Sym$-module of $\MR$ as
follows.
The basis element $\tilde S^I$ of $\BSym$, where $I=(i_0,i_1,\dots,i_r)$ is a
$B$-composition (that is, $i_0$ may be $0$), can be embedded as
\begin{equation}
\tilde S^I = S_{i_0}(A)S^{i_1i_2\dots i_r}(A|\Ab)\,.
\end{equation}
In the sequel, we identify $\BSym$ with its image under this embedding.

As in $\Sym$, one can define by triangularity the analog of the ribbon basis
(\cite{Chow}):
\begin{equation}
\tilde S^I = \sum_{J\leq I} \tilde R_J,
\end{equation}
where $J\leq I$ if the $B$-descent set of $J$ is a subset of the $B$-descent
set of $I$.
Note that we have in particular
$\tilde S^{0n} = \tilde R_{0n} + \tilde R_n$.

Note also that, thanks to that definition,
${S^I}^\#=\tilde S^I$ and, thanks to the transitions between all bases,
\begin{equation}
\label{ridiese}
R_I^\# = \tilde R_{0I} + \tilde R_I.
\end{equation}

\subsection{Type $B$ permutations and descents in $B_n$}

The hyperoctahedral group $B_n$ is the group of signed permutations. A signed
permutation can be denoted by $w=(\sigma,\epsilon)$ where $\sigma$ is an
ordinary permutation and $\epsilon\in\{\pm 1\}^n$, such that
$w(i)=\epsilon_i\sigma(i)$.
If we set $w(0)=0$, then, $i\in[0,n-1]$ is a $B$-descent of $w$ if
$w(i)>w(i+1)$.
Hence, the $B$-descent set of $w$ is a subset
$D=\{i_0,i_0+i_1,\dots,i_0+\dots+i_{r-1}\}$ of $[0,n-1]$.
We then associate with $D$ the type-$B$ composition 
$(i_0-0,i_1,\dots,i_{r-1},n-i_{r-1})$.

\subsection{Free quasi-symmetric functions of level 2}

Let us now move to generalizations of $\FQSym$. As in the case of $\NCSF$, the
most natural way is to change the usual alphabet into two alphabets, one of
positive letters and one of negative letters and to define a basis indexed by
signed permutations as a realization on words on both alphabets.
%
This algebra is $\FQSym^{(2)}$, the algebra of free quasi-symmetric functions
of level $2$, as defined in~\cite{NT}.

Let us set
\begin{gather}
A^{(0)}=A = \{a_1<a_2<\dots<a_n<\dots\}\,,\\
A^{(1)}=\Ab = \{\dots <\bar a_n <\dots <\bar a_2<\bar a_1\}\,,
\end{gather}
and order ${\bf A}=\bar A\cup A$ by $\bar a_i<a_j$ for all $i,j$. Let us
also denote by $\std$ the standardization of signed words with respect to this
order.

We shall also need the \emph{signed standardization} $\cstd$, defined as
follows.
Represent a signed word ${\bf w}\in{\bf A}^n$ by a pair $(w,\epsilon)$, where
$w\in A^n$ is the underlying unsigned word, and $\epsilon\in\{\pm 1\}^n$ is
the vector of signs. Then $\cstd(w,\epsilon)=(\std(w),\epsilon)$.

Then, $\FQSym^{(2)}$ is spanned by the polynomials in $A\cup\Ab$
\begin{equation}
\label{defG}
\G_{\sigma,u} := \sum_{w\in \A^n ; \cstd(w)=(\sigma,u)} w \quad\in\Z\free{\A}.
\end{equation}

Let $(\sigma',u')$ and $(\sigma'',u'')$ be signed permutations.
Then (see~\cite{NT,NT-super})
\begin{equation}
\label{multGb}
\G_{\sigma',u'}\,\,\G_{\sigma'',u''}
= \sum_{\sigma\in \sigma'\convol\sigma''} \G_{\sigma,u'\cdot u''}.
\end{equation}
We denote by $\moinsu$ the number of entries $-1$ in $\epsilon$.

\subsection{Free super-quasi-symmetric functions}

The second algebra generalizing the algebra $\FQSym$ is $\FQSym(A|\Ab)$. It
comes from the transformation of alphabets applied to $\FQSym$ as
$\NCSF(A|\Ab)$ comes from $\NCSF$. To do this, we first need to recall that
$\FQSym^{(2)}$ is equipped with an internal product.

Indeed, viewing signed permutations as elements of the group
$\{\pm 1\}\wr\SG_n$, we have the internal product
\begin{equation}
\G_{\alpha,\epsilon}* \G_{\beta,\eta}=\G_{(\beta,\eta)\circ(\alpha,\epsilon)}
=\G_{\beta\circ\alpha, (\eta\alpha)\cdot \epsilon},
\end{equation}
with $\eta\alpha=(\eta_{\alpha(1)},\dots,\eta_{\alpha(n)})$
and $\epsilon\cdot \eta=(\epsilon_1\eta_1,\dots,\epsilon_n\eta_n)$.

\smallskip
We can now embed $\FQSym$ into $\FQSym^{(2)}$ by
\begin{equation}
\G_\sigma \mapsto \G_{(\sigma,1^n)},
\end{equation}
which allows us to define
\begin{equation}
\G_\sigma^\# := \G_\sigma(A|\Ab)=\G_\sigma * \sigma_1^\#,
\end{equation}
so that $\FQSym(A|\Ab)$ is the algebra spanned by the $\G_\sigma(A|\Ab)$.

\begin{theorem}[\cite{NT-super}, Thm. 3.1]
\label{GAAb}
The expansion of $\G_\sigma(A|\Ab)$ on the basis $\G_{\tau,\epsilon}$ is
\begin{equation}\label{Gdiese}
\G_\sigma(A|\Ab)=
\sum_{\std(\tau,\epsilon)=\sigma}\G_{\tau,\epsilon}\,.
\end{equation}
\end{theorem}

\subsubsection{Embedding $\NCSF^\#$ and $\BSym$ into $\FQSym^{(2)}$}

One can embed $\BSym$ into $\FQSym^{(2)}$ as one embeds $\Sym$ into $\FQSym$
(see~\cite{Chow}) by
\begin{equation}
\label{tildeR}
\tilde R_{I} = \sum_{\Bdes(\pi)=I} \G_{\pi},
\end{equation}
where $I$ is any $B$-composition.

Given Equation~\eqref{ridiese} relating $R_I^\#$ and the $\tilde R_I$,
one has
\begin{equation}
\label{Rdiese}
R_I^\sharp = \sum_{\Des(\pi)=I} \G_{\pi},
\end{equation}
where $I$ is any (usual) composition.

\section{Algebraic theory in type $B$}
\label{section4}

\subsection{Alternating permutations of type $B$}

\subsubsection{Alternating shapes}

Let us say that a signed permutation $\pi\in B_n$ is \emph{alternating}
if $\pi_1\!<\!\pi_2\!>\!\pi_3\!<\!\dots$ (shape $2^m$ or $2^m1$).

\medskip
{\scriptsize
Here are the alternating permutations of type $B$ for $n\leq4$:
\begin{equation}
\bar 1, 1,
\quad
12, \bar12, \bar21, \bar2\bar1
\end{equation}
\begin{equation}
\label{alt-3}
12\bar3, \bar12\bar3, 132, 13\bar2, \bar132, \bar13\bar2, \bar21\bar3,
\bar2\bar1\bar3, 231, 23\bar1, \bar231, \bar23\bar1, \bar31\bar2,
\bar3\bar1\bar2, \bar321, \bar32\bar1,
\end{equation}
\smallskip
\begin{equation}
\begin{split}
&
12\bar34, \bar12\bar34, 12\bar43, 12\bar4\bar3, \bar12\bar43,
\bar12\bar4\bar3,
1324, 13\bar24, \bar1324, \bar13\bar24, 13\bar42, 13\bar4\bar2,
\\ &
\bar13\bar42, \bar13\bar4\bar2, 1423, 14\bar23, \bar1423,
\bar14\bar23,
14\bar32, 14\bar3\bar2, \bar14\bar32, \bar14\bar3\bar2,
\bar21\bar34, \bar2\bar1\bar34,
\\ &
\bar21\bar43, \bar21\bar4\bar3, \bar2\bar1\bar43,
\bar2\bar1\bar4\bar3, 2314, 23\bar14,
\bar2314, \bar23\bar14, 23\bar41, 23\bar4\bar1, \bar23\bar41,
\bar23\bar4\bar1,
\\ &
2413, 24\bar13, \bar2413,
\bar24\bar13, 24\bar31, 24\bar3\bar1,
\bar24\bar31, \bar24\bar3\bar1, \bar31\bar24,
\bar3\bar1\bar24, \bar31\bar42, \bar31\bar4\bar2,
\\ &
\bar3\bar1\bar42, \bar3\bar1\bar4\bar2, \bar3214, \bar32\bar14,
\bar32\bar41, \bar32\bar4\bar1,
\bar3\bar2\bar41, \bar3\bar2\bar4\bar1, 3412, 34\bar12, \bar3412,
\bar34\bar12,
\\ &
34\bar21, 34\bar2\bar1, \bar34\bar21,
\bar34\bar2\bar1,
\bar41\bar23, \bar4\bar1\bar23,
\bar41\bar32, \bar41\bar3\bar2, \bar4\bar1\bar32,
\bar4\bar1\bar3\bar2, \bar4213, \bar42\bar13,
\\ &
\bar42\bar31, \bar42\bar3\bar1, \bar4\bar2\bar31,
\bar4\bar2\bar3\bar1, \bar4312, \bar43\bar12,
\bar43\bar21, \bar43\bar2\bar1
\end{split}
\end{equation}
}

Hence, $\pi$ is alternating iff $\bar \pi$ is a
$\beta$-snake in the sense of \cite{Arnold}.
Hence, the sum in $\FQSym^{(2)}$ of all $\G_\pi$ labeled by
alternating signed permutations is, as already proved in~\cite{Chow}
\begin{equation}
\label{fqs2}
\X = (X+Y)^\# = \SEC^\#+\TAN^\#
   = \SEC^\# ( 1+\SIN^\#)
   =\sum_{m\ge 0}(R_{(2^m)}^\#+R_{(2^m1)}^\#).
\end{equation}

\subsubsection{Quasi-differential equations}

Let $d$ be the linear map acting on $\G_\pi$ as follows:
\begin{equation}
\label{dgb}
d\G_\pi=
\left\{
\begin{array}{cc}
\G_{uv}      & \text{if } \pi=unv, \\
\G_{u\bar v} & \text{if } \pi=u\bar nv.
\end{array}
\right.
\end{equation}
This map lifts to $\FQSym^{(2)}$ the derivation $\partial$
of \eqref{part}, although it is not itself a derivation.
We then have

\begin{theorem}
The series $\X$ satisfies the quasi-differential equation
\begin{equation}
\label{dx1p}
d\X = 1+ \X^2\,.
\end{equation}
\end{theorem}

\Proof
Indeed, let us compute what happens when applying $d$ to $R_{(2^m)}^\#$, the
property being the same with $d R_{(2^m1)}^\#$.
Let us fix a permutation $\sigma$ of shape $(2^m)$. If $n$ appears in
$\sigma$, let us write $\sigma = unv$.
Then $dG_\sigma$ appears in the product $\G_{\Std(u)}\G_{\Std(v)}$ and $u$ and
$v$ are of respective shapes $(2^n1)$ and $(2^{m-n-1})$. If $\overline n$
appears in $\sigma$, let us write again $\sigma = u\overline nv$.
Then $dG_\sigma$ appears in the product $\G_{\Std(u)}\G_{\Std(\overline v)}$
and $u$ and $\overline v$ are of respective shapes $(2^n)$ and $(2^{m-n-1}1)$.
Conversely, any permutation belonging to $u\convol v$ with $u$ and $v$ of
shapes $(2^n1)$ and $(2^{m-n-1})$ has a shape $(2^m)$ if one adds $n$ in
position $2n+2$. The same holds for the other product, hence proving the
statement.
\qed

This is not enough to characterize $\X$ but we have the analog of
fixed point equation~\eqref{eqbinX}
\begin{equation}
\label{fixB}
\X = 1 + \G_1 + B(\X,\X),
\end{equation}
where
\begin{equation}
B(\G_\alpha,\G_\beta)=
\begin{cases}
\sum_{\gamma=u(n+1)v,\ \cstd(u)=\alpha,\ \cstd(v)=\beta}\G_\gamma\
 &\text{if $|\alpha|$ is odd},\\
\sum_{\gamma=u(\overline{n+1})\bar v,\ \cstd(u)=\alpha,\ \cstd(v)=\beta}
 \G_\gamma\ &\text{if $|\alpha|$ is even}.\\
\end{cases}
\end{equation}
Indeed, applying $d$ to the fixed point equation brings back
Equation~\eqref{dx1p} and it is clear from the definition of $B$ that all
terms in $\B(\G_\alpha,\G_\beta)$ are alternating signed permutations.

Solving this equation by iterations gives back the results
of~\cite[Section 4]{JV}.
Indeed, the iteration of Equation~(\ref{fixB}) yields the solution
\begin{equation}
\label{X2cbt}
\X = \sum_{T\in \CBT} B_T(\G_0=1,\G_1),
\end{equation}
where, for a tree $T$, $B_T(a,b)$ is the result of the evaluation of all
expressions formed by labeling by $a$ or $b$ the leaves of $T$ and by $B$ its
internal nodes. This is indeed the same as the polynomials $P_n$ defined
in~\cite[Section 4]{JV} since one can interpret the $\G_0$ leaves as empty
leaves in this setting, the remaining nodes then corresponding to all
increasing trees of the same shape, as can be seen on the definition of the
operator~$B$.

\subsubsection{Alternating signed permutations counted by number of signs}

Under the specialization $\bar A=tA$, $\X$ goes to the series
\begin{equation}
X(t;A)=\sum_{I}\left( 
\sum_{\pi \ \text{alternating},\ C(\std(\pi))=I}t^{m(\pi)}
\right)R_I(A)
\end{equation}
where $m(\pi)$ is the number of negative letters of $\pi$.
If we further set $A=z\E$, we obtain
\begin{equation} \label{tanalogue}
x(t,z) = \frac{1+\sin((1+t)z)}{\cos((1+t)z)}
\end{equation}
which reduces to
\begin{equation}
x(1,z)
 = \frac{1+ \sin 2z}{\cos 2z} 
 = \frac{\cos z+\sin z}{\cos z-\sin z} 
\end{equation}
for $t=1$, thus giving a $t$-analogue different from the one of~\cite{JV}.

\subsubsection{A simple bijection}

From~\eqref{tanalogue}, we have
\begin{equation}
 \sum_{n} z^n  \sum_{\pi \text{ alternating in $B_n$}} t^{m(\pi)}
 = \sec((1+t)z) + \tan((1+t)z).
\end{equation}
But another immediate interpretation of the series in the right-hand side is
\begin{equation}
 \sum_{n} z^n 
     \sum_{\pi \text{\ s.t. $|\pi|$ is alternating in $A_n$} } t^{m(\pi)}
 = \sec((1+t)z) + \tan((1+t)z).
\end{equation}
It is thus in order to give a bijection proving the equality of the generating
functions.
Let $\pi$ be an alternating signed permutation. We can associate with $\pi$
the pair $(\std(\pi), \epsilon)$ where $\epsilon$ is the sign vector such
that $\epsilon_i=1$ if $\pi^{-1}(i)>0$ and $\epsilon_i=-1$ otherwise.
The image of $\{1,\dots,n\}$ by $\pi$ is $\{ \epsilon_i i : 1\leq i \leq n \}$.
Since $\pi$ can be recovered from $\std(\pi)$ and the image of
$\{1,\dots,n\}$, this map is a bijection between signed alternating
permutations and pairs $(\sigma,\epsilon)$ where $\sigma$ is alternating and
$\epsilon$ is a sign vector. Then, with such a pair $(\sigma,\epsilon)$,
one can associate a signed permutation $\tau$ such that $|\tau|$ is
alternating simply by taking $\tau_i = \sigma_i\epsilon_i$.
The composition $\pi\mapsto (\sigma,\epsilon) \mapsto \tau$ gives the desired
bijection.

For example, here follow the 16 permutations obtained by applying the
bijection to the 16 alternating permutations of size $3$
(see Equation~\eqref{alt-3}):
\begin{equation}
\label{alt-3b}
23\bar1, \bar23\bar1, 132, 2\bar31, \bar132, \bar2\bar31, 2\bar3\bar1,
\bar2\bar3\bar1, 231, \bar231, 1\bar32, \bar1\bar32, 1\bar3\bar2,
\bar1\bar3\bar2, 13\bar2, \bar13\bar2.
\end{equation}

\subsection{Type $B$ snakes}

\subsubsection{An alternative version}

The above considerations suggest a new definition of type $B$ snakes, which is
a slight variation of the definition of~\cite{Arnold}. We want to end up with
the generating series
\begin{equation}
y(1,z) = \frac{1}{\cos z-\sin z} = \frac{\cos z + \sin z}{\cos 2z}
\end{equation}
after the same sequence of specializations. A natural choice, 
simple enough and given by a series in $\BSym$, is to set
\begin{equation}
\Y = (\COS+\SIN)\cdot \SEC^\#
= \left( \sum_{k\ge 0}(-1)^k(S_{2k}+S_{2k+1}) \right)
  \cdot \sum_{n\ge 0} R_{2^n}^\# .
\end{equation}
Now, $\Y$ lives in $\BSym$ and expands in the ribbon basis $\tilde R$ of
$\BSym$ as
\begin{equation}
\begin{split}
\Y
&= \left( \sum_{k\ge 0}(-1)^k(\tilde R_{2k}+\tilde R_{2k+1}) \right)
  \sum_{n\ge 0} R_{2^n}^\# \\
&= \sum_{n\geq0} (R_{2^n}^\# + \tilde R_{12^n}+\tilde R_{32^{n}}) \\
&\ 
 + \sum_{k\geq1;n\geq0} (-1)^k
 ( \tilde R_{2k2^n}+\tilde R_{2k+22^{n-1}}
  +\tilde R_{2k+12^n}+\tilde R_{2k+32^{n-1}})\\
&= \sum_{n\geq0} (R_{2^n}^\# + \tilde R_{12^n}+\tilde R_{32^{n}})
 - \sum_{n\geq0} (\tilde R_{2^{n+1}} + \tilde R_{32^n})
\end{split}
\end{equation}
which simplifies into
\begin{equation}
\label{Y1b}
\Y = 1+\sum_{n\ge 0}(\tilde R_{12^n} +\tilde R_{02^{n+1}}).
\end{equation}

In $\FQSym^{(2)}$, this is the sum of all $\G_\pi$ such that
\begin{equation}
\left\{
\begin{array}{ll}
0>\pi_1<\pi_2>\dots & \text{if $n$ is even}, \\
0<\pi_1>\pi_2<\dots & \text{if $n$ is odd}.
\end{array}
\right.
\end{equation}
Thus, for $n$ odd, $\pi$ is exactly a $B_n$-snake in the sense
of~\cite{Arnold}, and for $n$ even, $\bar\pi$ is a $B_n$-snake.
Clearly, the number of sign changes or of minus signs in snakes and in
these modified snakes are related in a trivial way so we have generating
series for both statistics in all cases.

\medskip
{\scriptsize
Here are these modified snakes for $n\le 4$:
\begin{equation}
1, \qquad
\bar12, \bar21, \bar2\bar1,
\end{equation}
\begin{equation}
1\bar23, 1\bar32, 1\bar3\bar2, 213, 2\bar13, 2\bar31,
2\bar3\bar1, 312, 3\bar12, 3\bar21, 3\bar2\bar1,
\end{equation}

\begin{equation}
\begin{split}
&
\bar12\bar34, \bar12\bar43, \bar12\bar4\bar3, \bar1324, \bar13\bar24,
 \bar13\bar42,
\bar13\bar4\bar2, \bar1423, \bar14\bar23, \bar14\bar32, \bar14\bar3\bar2,
 \bar21\bar34,
\\
&
   \bar2\bar1\bar34, \bar21\bar43, \bar21\bar4\bar3,
\bar2\bar1\bar43, \bar2\bar1
\bar4\bar3, \bar2314,
   \bar23\bar14, \bar23\bar41, \bar23\bar4\bar1,
\bar2413, \bar24\bar1
3, \bar24\bar31,
\\
&
   \bar24\bar3\bar1, \bar31\bar24, \bar3\bar1\bar24,
\bar31\bar42, \bar31
\bar4\bar2, \bar3\bar1\bar42,
   \bar3\bar1\bar4\bar2, \bar3214, \bar32\bar14,
\bar32\bar41, \bar32
\bar4\bar1, \bar3\bar2\bar41,
\\
&
   \bar3\bar2\bar4\bar1, \bar3412, \bar34\bar12,
\bar34\bar21, \bar34
\bar2\bar1, \bar41\bar23,
   \bar4\bar1\bar23, \bar41\bar32, \bar41\bar3\bar2,
\bar4\bar1\bar32, \bar4\bar1
\bar3\bar2, \bar4213,
\\
&
   \bar42\bar13, \bar42\bar31, \bar42\bar3\bar1,
\bar4\bar2\bar31, \bar4\bar2
\bar3\bar1, \bar4312,
   \bar43\bar12, \bar43\bar21, \bar43\bar2\bar1.
\end{split}
\end{equation}
}

\subsubsection{Snakes as particular alternating permutations}
\label{snakesssens}

Note that in the previous definition, snakes are not alternating permutations
for odd $n$.
So, instead, let us consider the generating series 
\begin{equation}
\Y = \COS\cdot \SEC^\# + \overline\SIN\cdot{\overline\SEC^\#}\,,
\end{equation}
where $f\mapsto\bar f$ is the involution of $\FQSym^{(2)}$ inverting the signs
of permutations.

Expanding $\Y$ in the $\tilde R$ basis, one gets
\begin{equation}
\label{bsnakes}
\Y = 1+\sum_{n\ge 0}(\tilde R_{02^n1} +\tilde R_{02^{n+1}}).
\end{equation}

As for type $B$ alternating permutations (see Equation~\eqref{dx1p}), the
series $Y$ satisfies a differential equation with the same linear map $d$ as
before (see Equation~\eqref{dgb}):

\begin{equation}
\label{dy1p}
d \Y = \Y\, \X.
\end{equation}

It is then easy to see that $\Y$ also satisfies
a fixed point equation similar to~\eqref{eqbinX}:
\begin{equation}
\label{fixBY}
\Y = 1 + B(\Y,\X).
\end{equation}

The iteration of~(\ref{fixBY})  brings up a solution close to~(\ref{X2cbt}):
\begin{equation}
\label{Y2cbt}
\Y = \sum_{T\in \CBT} B_T(\G_0=1,\G_1),
\end{equation}
where, for a tree $T$, $B_T(a,b)$ is now the result of the evaluation of all
expressions formed by labeling by $a$ or $b$ the leaves of $T$ and by $B$ its 
internal nodes. Note that in this case, the first leaf needs to have label
$a$.
This is the same as the trees defined in~\cite[Section 4]{JV}
since one can again interpret the $\G_0$ leaves as empty leaves in this
setting, the remaining nodes then corresponding to all increasing trees of the
same shape.

\subsubsection{Snakes from~\cite{Arnold}}

The generating series of the snakes of \cite{Arnold}, also in~$\BSym$ is
\begin{equation}
\overline{\COS}\cdot\overline{\SEC^\#}+\SIN\cdot\SEC^\#\,,
\end{equation}
and can be written as
\begin{equation}
\Y = 1+\sum_{n\ge 0}(\tilde R_{12^n} +\tilde R_{12^n1})
\end{equation}
on the ribbon basis.

The lift of the differential equation for $y(1)$ is given by a map $\delta$
similar to $d$, with $\delta unv=\bar u v$ and $\delta u\bar n v=u\bar v$.
Then
\begin{equation}
\delta \Y =\Y\X
\end{equation}
and we have a fixed point equation
\begin{equation}
\label{fixBs}
\Y = 1 + \hat B(\Y,\X)
\end{equation}
for an appropriate bilinear map  $\hat B$.

\section{Another combinatorial model}
\label{section5}

\subsection{An analogue of $\cos z-\sin z$ in $\FQSym^{(2)}$}

\begin{definition}
A signed permutation $\pi\in B_n$ is a \emph{valley-signed permutation} if,
for any $i\in[n]$, $\pi(i)<0$ implies that
\begin{itemize}
\item either $i>2$, and $\pi_{i-1}>0$, and $|\pi_{i-2}|>\pi_{i-1}<| \pi_i|$,
\item or $i=2$,  and $0<\pi_1<|\pi_2|$.
\end{itemize}
We denote by $\mathcal{V}_n$ the set of valley-signed permutations of size $n$.
\end{definition}

{\scriptsize
Here are these signed permutations, up to $n=4$:
\begin{equation}
1, \qquad
12, 1\bar2,  21,
\end{equation}
\smallskip
\begin{equation}
123, 1\bar23, 132, 1\bar32, 213, 21\bar3, 231, 2\bar31, 312, 31\bar2, 321,
\end{equation}
\smallskip
\begin{equation}
\begin{split}
&
1234, 1\bar234, 1243, 1\bar243, 1324, 132\bar4, 1\bar324, 1\bar32\bar4, 1342,
1\bar342, 1423, 142\bar3,
\\ &
1\bar423, 1\bar42\bar3, 1432, 1\bar432, 2134, 21\bar34, 2143, 21\bar43, 2314,
231\bar4, 2\bar314, 2\bar31\bar4,
\\ &
2341, 2\bar341, 2413, 241\bar3, 2\bar413, 2\bar41\bar3, 2431, 2\bar431, 3124,
31\bar24, 3142, 31\bar42,
\\ &
3214, 321\bar4, 3241, 32\bar41, 3412, 341\bar2, 3\bar412, 3\bar41\bar2, 3421,
3\bar421, 4123, 41\bar23,
\\ &
4132, 41\bar32, 4213, 421\bar3, 4231, 42\bar31, 4312, 431\bar2, 4321
\end{split}
\end{equation}
}

The terminology is explained by the following remark. Let
$\sigma\in\mathfrak{S}_n$, and let us examine how to build a valley-signed
permutation $\pi$ so that $\pi_i=\pm \sigma_i$.
It turns out that for each valley
$\sigma(i-1)>\sigma(i)<\sigma(i+1)$, we can choose independently the sign of
$\pi(i+1)$ (here $1\leq i <n$ and $1$ is a valley if $\sigma(1)<\sigma(2)$).

The goal of this section is to obtain the noncommutative generating functions
for the sets $\mathcal{V}_n$.

\begin{theorem}
\label{um1}
The following series
\begin{equation}
U= 1-\left[
\G_1 + \G_{1\bar 2} - \G_{12\bar 3} - \G_{1\bar 2 3\bar 4}
+\G_{12\bar 3 4\bar 5}+\G_{1\bar 2 3\bar 4 5\bar6}+\dots
\right]
\end{equation}
is again a lift of $\cos z-\sin z$ in $\FQSym^{(2)}$. It satisfies
\begin{equation}
U^{-1}
= \sum_{n\geq0}  \sum_{\pi \in \mathcal{V}_n } \G_\pi.
\end{equation}
Hence the $\pi$ occuring in this expansion are in bijection with snakes of
type $B$.
\end{theorem}

This result is a consequence of the next two propositions.

\begin{definition}
Let $\mathcal{R}_{2n}\subset B_{2n}$ be the set of signed permutations $\pi$
of size $2n$ such that $ |\pi|$ is of shape $2^n$,
and for any $1\leq i \leq n$, we have $\pi(i)>0$ iff $i$ is odd.
Let
\begin{equation}
V = \sum_{n\geq 0} \sum_{\pi\in\mathcal{R}_{2n}} \G_\pi
  =  \G_\epsilon + \G_{1\bar2} + \G_{1\bar32\bar4} +  \G_{1\bar42\bar3}+\dots
\end{equation}
Let $\mathcal{R}_{2n+1}\subset B_{2n+1}$ be the set of signed permutations
$\pi$ of size $2n+1$ such that $ |\pi| $ is of shape $12^n$, $\pi_1>0$, and
for any $2\leq i \leq n$, we have $\pi(i)>0$ iff $i$ is even.
Let
\begin{equation}
W = \sum_{n\geq 0} \sum_{\pi\in\mathcal{R}_{2n+1}} \G_\pi
  =  \G_1  + \G_{21\bar 3} + \G_{31\bar 2} + \G_{21\bar43\bar5} + \dots
\end{equation}
\end{definition}

Note that $\mathcal{R}_n\subset\mathcal{V}_n$.
Clearly, $\# \mathcal{R}_n$ is the number of alternating permutations of
$\mathfrak{S}_n$ since, given $|\pi|$, there is only one possible choice for
the signs of each $\pi_i$.
So $V$ and $W$ respectively lift $\sec$ and $\tan$ in $\FQSym^{(2)}$.
Now, given that the product rule of the $\G_\sigma$ does not affect the signs,
a simple adaptation of the proof of the $A_n$ case shows:

\begin{proposition}
We have
\begin{align}
 V^{-1}   &=  1  - \G_{1\bar2} + \G_{1\bar23\bar4} - \G_{1\bar23\bar45\bar6}
               + \dots, \\
 WV^{-1}  &= \G_1 -  \G_{12\bar3} + \G_{12 \bar34\bar5} -
             \G_{12\bar34\bar56\bar7} + \dots
\end{align}
\end{proposition}
\hfill\qed

Note that $U=(1-W)V^{-1}$. So, to complete the proof of the theorem,
it remains to show:

\begin{proposition}
We have
\begin{equation}
  V(1-W)^{-1} =  \sum_{n\geq0}  \sum_{\pi \in \mathcal{V}_n } \G_\pi.
\end{equation}
\end{proposition}

\Proof
We can write $V(1-W)^{-1} = V+VW+VW^2+\dots$, and expand everything in terms
of the $\G_\pi$, using their product rule (see Equation~\eqref{multGb}). We 
obtain a sum of $\G_{u_1 \dots u_k}$ where $\cstd(u_1) \in \mathcal{R}_{2*}$
and $\cstd(u_i)\in \mathcal{R}_{2*+1}$ for any $i\geq2$ (where $*$ is any
integer).
The sum is \emph{a priori} over lists $(u_1,\dots,u_k)$ such that
$u=u_1\dots u_k\in B_n$.
Actually, if the words $u_1,\dots,u_k$ satisfy the previous conditions, then
$u$ is in $\mathcal{V}_{n}$. Indeed, the first two letters of each $u_i$ are
not signed, each $u_i$ is a valley-signed permutation, which implies that
$u$ itself is a valley-signed permutation.

Conversely, it remains to show that this factorization exists and is unique
for any $u\in\mathcal{V}_{n}$. First, observe that in $u=u_1 \dots u_k$,
$u_k$ is the only suffix of odd length having the same signs as an element
of $\mathcal{R}_{2*+1}$ (since the first two letters are positive and the
other alternate in signs, it cannot be itself a strict suffix of an element
of $\mathcal{R}_{2*+1}$). So the factorization can be obtained by scanning $u$
from right to left.
\qed

Let 
\begin{equation}
Z := \G_{\bar1} - \G_{\bar12} - \G_{\bar12\bar3} + \G_{\bar123\bar4}
    +\G_{\bar12\bar 3 4\bar 5} - \G_{\bar123\bar45\bar6} + \dots
\end{equation}
Precisely, the $n$-th term in this expansion is the permutation $\sigma$
of $B_n$ such that, $|\sigma|=id$, $\sigma(1)=\bar1$, $\sigma(2)=2$ and, for
$3\leq i\leq n$, $\sigma(i)<0$ iff $n-i$ is even. The sign corresponds to the
sign of $z^n$ in the expansion of $\cos z+\sin z-1$, so that
it is a lift of $\cos z+\sin z-1$ in $\FQSym^{(2)}$.

\begin{theorem}
\label{56}
The series $Z\,U^{-1}$ is a sum of $\G_\pi$ without multiplicities in
the Hopf algebra $\FQSym^{(2)}$. Hence the $\pi$ occuring in this expansion
are in bijection with snakes of type $D$ (see Section~\ref{sectionD}).

Moreover, the series $(1+Z)\,U^{-1}$ is also a sum of $\G_\pi$ without
multiplicities in the Hopf algebra $\FQSym^{(2)}$. Hence the $\pi$ occuring in
this expansion are in bijection with alternating permutations of type $B$.
\end{theorem}

\medskip
{\scriptsize
Here are the elements of $Z\,U^{-1}$, up to $n=4$:

\begin{equation}
\bar1, \qquad
\bar21,
\end{equation}
\begin{equation}
\bar213, \bar21\bar3, \bar312, \bar31\bar2, \bar321,
\end{equation}
\begin{equation}
\begin{split}
&
\bar21\bar34, \bar2134, \bar21\bar43, \bar2143, \bar31\bar24, \bar3124,
\bar31\bar42, \bar3142, \bar321\bar4, \bar3214, \bar32\bar41, \bar3241,
\\
&
 \bar41\bar23, \bar4123, \bar41\bar32, \bar4132, \bar421\bar3, \bar4213,
 \bar42\bar31, \bar4231, \bar431\bar2, \bar4312, \bar4321.
\end{split}
\end{equation}
}

\medskip
\Proof
Let $\varepsilon_i$ be the linear operator sending a word $w$ to the word
$w'$ where $w'$ is obtained from $w$ by sending its $i$th letter to its
opposite and not changing the other letters.

Then if one writes $1+Z=E+O$ as a sum of an even and an odd series, one has
\begin{equation}
E=\varepsilon_1\varepsilon_2 V^{-1} \qquad O=\epsilon_1(W V^{-1}).
\end{equation}
Since $\varepsilon_i(ST)=\varepsilon_i(S)T$ if $S$ contains only terms of
size at least $i$, an easy rewriting shows that
\begin{equation}
(1+Z)V = V + \varepsilon_1(W) + \varepsilon_1\varepsilon_2(1-V).
\end{equation}
Then, one gets
\begin{equation}
(1+Z)V(1-W)^{-1} = U^{-1} + \varepsilon_1(W+\varepsilon_2(1-V))(1-W)^{-1}.
\end{equation}
So it only remains to prove that
$Q=\varepsilon_1(W+\varepsilon_2(1-V))(1-W)^{-1}$ has only positive terms and
has no term in common with $U^{-1}$.
This last fact follows from the fact that all terms in $Q$ have a negative
number as first value.
Let us now prove that $Q$ has only positive terms. First note that
$W+\varepsilon_2(1-V)$ is an alternating sum of permutations of shapes $2^n$
and $1 2^n$. Hence, any permutation of shape $2^n$ can be associated with a
permutation of shape $12^{n-1}$ by removing its first entry and standardizing
the corresponding word.
Now, 
all negative terms $-\G_v \G_w$ come from $-\varepsilon_2(V)(1-W)^{-1}$ and
are annihilated by the term $\G_{v'}\G_1\G_w$ where $v'$ is obtained from $v$
by the removal-and-standardization process described before.
\qed

\subsection{Another proof of Theorem~\ref{um1}}

Let $\Sg_n$ be the group algebra of $\{\pm 1\}^n$. We identify
a tuple of signs with a word in the two symbols $1,\bar 1$, and
the direct sum of the $\Sg_n$ with the free associative algebra
on these symbols.

We can now define the algebra of
\emph{signed noncommutative symmetric functions} as
\begin{equation}
\sSym := \bigoplus_{n\ge 0} \Sym_n\otimes \Sg_n
\end{equation}
endowed with the product
\begin{equation}
(f\otimes u)\cdot (g\otimes v)=fg\otimes uv\,.
\end{equation}
It is naturally embedded in $\FQSym^{(2)}$ by
\begin{equation}
R_I\otimes u =\sum_{C(\sigma)=I}\G_{\sigma,u}\,.
\end{equation}
With this at hand, writing $U=P-Q$ with
\begin{equation}
P=\sum_{m\ge 0}(-1)^mR_{2m}\otimes (1\bar1)^m
\quad\text{and}\quad
Q=\sum_{m\ge 0}(-1)^mR_{2m+1}\otimes 1(1\bar 1)^m,
\end{equation}
it is clear that
\begin{equation}
P^{-1}=\sum_{m\ge 0}R_{(2^m)}\otimes(1\bar 1)^m = V
\ \text{and}\
W=\sum_{m\ge 0}R_{(12^m)}\otimes 1(1\bar 1)^m = QV,
\end{equation}
so that $U=(1-W)V^{-1}$, and $U^{-1}=V(1-W)^{-1}$ can now be computed by
observing that
\begin{equation}
(1-W)^{-1}=\sum_I R_I\otimes p_I(1,\bar 1)
\end{equation}
where $p_I$ is the sum of words in $1,\bar 1$ obtained by associating with
each tiling of the ribbon shape of $I$ by tiles of shapes
$(1,2^{m_1},\dots,1,2^{m_k})$
the word $1(1\bar 1)^{m_1}\dots 1(1\bar 1)^{m_k}$.
For example, $p_{12}=111+11\bar 1$,
for there are two tilings, one by three shapes $(1)$ and one by $(12)$.
The proof of Theorem~\ref{56} can be reformulated in the same way.

\subsection{Some new identities in $\BSym$}

Let us have a closer look at the map $\bar A\mapsto tA$.

By definition, $R_I(A|tA)=R_I((1-q)A)|_{q=-t}$, so that the generating
series for one-part compositions is
\begin{equation}
\label{super}
\sum_{n\ge 0}R_n((1-q)A)x^n = \sigma_x((1-q)A)=\lambda_{-qx}(A)\sigma_x(A)\,.
\end{equation}
One can now expand this formula in different bases.
Tables are given in Section~\ref{ann-tilde}.

\subsubsection{Image of the $\tilde R$ in the $S$ basis of $\NCSF$}

The first identity gives the image of a ribbon in the $S$ basis:
\begin{equation}
\label{newid1}
R_I(A|tA)
= \sum_{J} (-1)^{l(I)+l(J)}\ (1-(-t)^{j_r})\ (-t)^{\sum_{A(I,J)} j_k}
  S^J
\end{equation}
where $A(I,J)=\{p| j_1+\dots+j_p \not\in \Des(I) \}$.
Indeed, expanding $\lambda_{-qx}(A)$ on the basis $S^J$, one finds
\begin{equation}
\sum_{n\ge 0}R_n((1-q)A)x^n
=\sum_{n\ge 0}x^n
      \left(\sum_{k+m=n}(-q)^k\sum_{K\vDash k}(-1)^{k-\ell(K)}S^KS_m \right)
\end{equation}
Each composition of a given $n$ occurs twice in the sum, so that
\begin{equation}
R_n((1-q)A)=\sum_{J\vDash n}(-1)^{1-\ell(J)}(q^{n-j_r}-q^n)S^J\,.
\end{equation}
Hence, Equation~(\ref{newid1}) is true for one-part compositions. The general
case follows by induction from the product formula
\begin{equation}
R_I R_j = R_{Ij} + R_{I\triangleright j}\,.
\end{equation}

Now, as for the $\tilde{R}$, we have

\begin{proposition}
Let $I$ be a type-$B$ composition and set $I=:(i_0,I')$.
Then the expansion of the
$\tilde{R}_I(A,tA) $
on the $S^J(A)$ is
\begin{equation}
\left\{
\begin{array}{ll}
\displaystyle
(-1)^{\ell(I')+n} t^nS^n +
\sum_{J\in \II_1}
(-1)^{\ell(I)+\ell(J)}\ (1-(-t)^{j_r})\ (-t)^{\sum_{A(I,J)} j_k} S^J
& \text{if $i_0=0$,} \\
\displaystyle
(-1)^{\ell(I)+1} S^n +
\sum_{J\in \II_2}
(-1)^{\ell(I)+\ell(J)}\ (1-(-t)^{j_r})\ (-t)^{\sum_{A(I,J)} j_k} S^J
& \text{otherwise},
\end{array}
\right.
\end{equation}
where $\II_2$ is the set of compositions of $n$ different from $(n)$ whose
first part is a sum of any prefix of $I$,
and $\II_1$ the set complementary to $\II_2$ in all compositions of $n$
different from $(n)$.
\end{proposition}

For example, with $I=(1321)$, $\II_2(I)$ is the set of compositions of $7$
whose first part is either $1$, $4$, or $6$.

\medskip
\Proof
By induction on the length of $I$.
First, $R_n(A,tA)=R_n(A)=S_n(A)$ and
\begin{equation}
\tilde R_n(A,tA)=R_n(A) \text{\quad and\quad}
\tilde R_{0n}(A,tA) = R_n(A|tA) - R_n(A).
\end{equation}
Now, the formula
$\tilde R_I R_j =\tilde R_{Ij} + \tilde R_{I\triangleright j}$
together with Equation~(\ref{newid1}) implies the general case.
\qed

Note that this also means that one can compute the matrices recursively.
Indeed, if one denotes by $K_0(n)$ (resp. $K_1(n)$) the matrix expanding the
$\tilde R_{I}(A,tA)$ where $i_0=0$ (resp. $i_0\not=0$) on the $S^J$, one has the
following structure:
\begin{equation}
\left( K_0(n\!+\!1)\, K_1(n\!+\!1) \right) =
\left(
\begin{array}{cccc}
-tK_0(n)             & tK_0(n) & K_1(n) & -K_1(n) \\
t(K_0(n)\!+\!K_1(n)) & 0       & 0      & K_0(n)\!+\!K_1(n)
\end{array}
\right)
\end{equation}

\subsubsection{Image of the $\tilde R$ in the $\Lambda$ basis of $\NCSF$}

Expanding Equation~(\ref{super}) on the basis $\Lambda^J$, the same reasoning
gives as well
\begin{equation}
\label{newid2}
R_I(A|tA)
= \sum_{J} (-1)^{n+1+l(I)+l(J)}\ (1-(-t)^{j_1})\
  (-t)^{\sum_{k\in A'(I,J)} j_k}
  \Lambda^J
\end{equation}
where $A'(I,J)=\{\ell| j_1+\dots+j_{\ell-1} \in \Des(I) \}$.

Now, as for the  $\tilde{R}$, we have

\begin{proposition}
Let $I$ be a type-$B$ composition and set $I=:(i_0,I')$.
Then the expansion of the
$\tilde{R}_I(A,tA) $
on the $\Lambda^J(A)$ is
\begin{equation}
\left\{
\begin{array}{cc}
\displaystyle
\sum_{J} (-1)^{n+l(I')+l(J)}\ (-t)^{j_1+\sum_{k\in A'(I',J)} j_k} \Lambda^J
& \text{if $i_0=0$,} \\
\displaystyle
\sum_{J} (-1)^{n+1+l(I')+l(J)}\ (-t)^{\sum_{k\in A'(I',J)} j_k} \Lambda^J
& \text{otherwise}.
\end{array}
\right.
\end{equation}
\end{proposition}
Note that this is coherent with the fact that $0\in\Des(I)$ if $i_0=0$.

\medskip
\Proof
The proof is again by induction of the length of $I$ following the same steps
as in the expansion on the $S$.
\qed

Again, one can compute the matrices recursively.
If one denotes by $L_0(n)$ (resp. $L_1(n)$) the matrix expanding the
$\tilde R_{I}(A,tA)$ where $i_0=0$ (resp. $i_0\not=0$) on the $\Lambda^J$, one
has the following structure:
\begin{equation}
\left( L_0(n\!+\!1)\, L_1(n\!+\!1) \right) =
\left(
\begin{array}{cccc}
tL_0(n) & -tL_0(n) & -L_1(n) & L_1(n) \\
tL_1(n) &  tL_0(n) &  L_1(n) & L_0(n)
\end{array}
\right).
\end{equation}

\subsubsection{Image of the $\tilde R$ in the $R$ basis of $\NCSF$}

The expansion of $R_I(A|tA)$ has been discussed in~\cite{NCSF2}
and~\cite{NT-super}. Recall that a \emph{peak} of a composition is a cell of
its ribbon diagram having no cell to its right nor on its top (compositions
with one part have by convention no peaks) and that a \emph{valley} is a cell
having no cell to its left nor at its bottom.

The formula is the following:
\begin{equation}
\label{RIsup}
R_I(A|tA) = \sum_{J} (1+t)^{v(J)} t^{b(I,J)} R_J(A),
\end{equation}
where the sum is over all compositions $J$ such that $I$ has either a peak or
a valley at each peak of $J$.
Here $v(J)$ is the number of valleys of $J$ and $b(I,J)$ is the number of
values $d$ such that, either $d$ is a descent of $J$ and not a descent of $I$,
or $d-1$ is a descent of $I$ and not a descent of $J$.

In the case of the $\tilde R$, the matrices satisfy a simple induction.
If one denotes by $M_0(n)$ (resp. $M_1(n)$) the matrix expanding the
$\tilde R_{I}(A,tA)$ where $i_0=0$ (resp. $i_0\not=0$) on the $R$, one has the
following structure which follows directly from the interpretation in terms of
signed permutations:
\begin{equation}
\left( M_0(n\!+\!1)\, M_1(n\!+\!1) \right) =
\left(
\begin{array}{cccc}
tM_1(n)          & tM_0(n) & M_1(n) & M_0(n) \\
t(M_0(n)\!+\!M_1(n)) & 0       & 0      & M_0(n)\!+\!M_1(n)
\end{array}
\right).
\end{equation}

For example, one can check this result on Figure~\ref{fig-R-23}.
One then recovers the matrix of $R_I(A|tA)$ on the $R$ as $M_0+M_1$.

\section{Euler-Bernoulli triangles}
\label{section6}

\subsection{Alternating permutations of type $B$}

Counting ordinary (type $A$) alternating permutations according to
 their last value yields
the Euler-Bernoulli triangle, sequence~A010094 or~A008281 of~\cite{Slo}
depending on whether one requires a rise or a descent at the first position.

The same can be done in type $B$ for alternating permutations and snakes.
Since usual snakes begin with a descent, we shall count type $B$ permutations
of ribbon shape $2^m$ or $2^m1$ according to their last value.
We then get the table

\bigskip
\begin{equation}
\label{AltBtriang}
\begin{array}{c|ccccccccccccc}
n\backslash p&-6 &-5 &-4 &-3 &-2 &-1 & 0 & 1 & 2 & 3 & 4 & 5 & 6  \\
\hline
 1 &   &   &   &   &   & 1 & 0 & 1 &   &   &   &   &    \\
 2 &   &   &   &   & 0 & 1 & 0 & 1 & 2 &   &   &   &    \\
 3 &   &   &   & 4 & 4 & 3 & 0 & 3 & 2 & 0 &   &   &    \\
 4 &   &   & 0 & 4 & 8 &11 & 0 &11 &14 &16 &16 &   &    \\
 5 &   &80 & 80&76 &68 &57 & 0 &57 &46 &32 &16 & 0 &    \\
 6 & 0 & 80&160&236&304&361& 0 &361&418&464&496&512&512 \\
\end{array}
\end{equation}

\bigskip
\begin{proposition}
\label{tri-alt}
The table counting type $B$ alternating permutations by their last value
is obtained by the following algorithm: first separate the picture by the
column $p=0$ and then compute two triangles. Put $1$ at the top of each
triangle and compute the rest as follows:
fill the second row of the left (resp. right) triangle as the sum of the
elements of the first row (resp. strictly) to their \emph{left}.
Then fill the third row of the right (resp. left) triangle as the sum of the
elements of the previous row (resp. strictly) to their \emph{right}. Compute
all rows successively by reading from left to right and right to left
alternatively.
\end{proposition}

This is the analogue for alternating permutations of Arnol'd's construction for
snakes of type $B$~\cite{Arnold}.

\medskip
\Proof
Let  $S(n,p)$ be the set of alternating permutations of  $B_n$
ending with~$p$.

The proof is almost exactly the same as for type $A$, with one exception: it
is obvious that $S(n,1)=S(n,-1)$.
Since the reading order changes from odd rows to even rows, let us assume 
that $n$ is even and consider both sets $S(n,p)$ and $S(n,p-1)$.
The natural injective map of $S(n,p-1)$ into $S(n,p)$ is simple: exchange
$p-1$ with $p$ while leaving the possible sign in place.
The elements of $S(n,p)$ that were not obtained previously are the
permutations ending by $p-1$ followed by $p$. Now, removing $p-1$ and
relabeling the remaining elements in order to get a type $B$ permutation, one
gets elements that are in bijection with elements of either $S(n-1,p)$ or
$S(n-1,p-1)$, depending on the sign of $p$.
\qed

\subsection{Snakes of type $B$}

The classical algorithm computing the number of type $B$ snakes (also known as
Springer numbers, see Sequence~A001586 of~\cite{Slo}) makes use of the
double Euler-Bernoulli triangle.

\begin{proposition}
\label{tri-sna}
The table counting snakes of type $B$ by their last value
is obtained by the following algorithm: first separate the picture by the
column $p=0$ and then compute two triangles. Put $1$ at the top of the left
triangle and $0$ at the top of the right one
and compute the rest as follows:
fill the second row of the left (resp. right) triangle as the sum of the
elements of the first row (resp. strictly) to their \emph{left}.
Then fill the third row of the right (resp. left) triangle as the sum of the
elements of the previous row (resp. strictly) to their \emph{right}. Compute
all rows successively by reading from left to right and right to left
alternatively.
\end{proposition}

Here are the first rows of both triangles:
\bigskip
\begin{equation}
\label{EBtriang}
\begin{array}{c|ccccccccccccc}
n\backslash p
   &-6 &-5 &-4 &-3 &-2 &-1 &   & 1 & 2 & 3 & 4 &5  &6 \\
\hline
 1 &   &   &   &   &   & 1 &   & 0 &   &   &   &   &  \\
 2 &   &   &   &   & 0 & 1 &   & 1 & 1 &   &   &   &  \\
 3 &   &   &   & 3 & 3 & 2 &   & 2 & 1 & 0 &   &   &  \\
 4 &   &   & 0 & 3 & 6 & 8 &   & 8 &10 &11 &11 &   &  \\
 5 &   &57 &57 &54 &48 &40 &   &40 &32 &22 &11 & 0 &  \\
 6 & 0 &57 &114&168&216&256&   &256&296&328&350&361&361\\
\end{array}
\end{equation}

\bigskip
\Proof
The proof is essentially the same as in the case of alternating permutations
of type $B$: it amounts to a bijection between a set of snakes on the one side
and two sets of snakes on the other side.
\qed

One also sees that each row of the triangles of the alternating type
$B$ permutations presented in Equation~(\ref{AltBtriang})
can be obtained, up to reversal, by adding or subtracting the mirror image of
the left triangle to the right triangle.
For example, on the fifth row, the sums are $40+40=80$, then $48+32=80$, then
$54+22=76$, then $57+11=68$, and $57+0=57$; the differences are
$57-0=57$, then $57-11=46$, then $54-22=32$, then $48-32=16$, and $40-40=0$.
These properties follow from the  induction patterns.

These numerical properties indicate that one can
split alternating permutations ending with $(-1)^{n-1}i$ into two sets and
obtain alternating permutations beginning with $(-1)^{n}i$ by somehow taking
the "difference" of these two sets.
On the alternating permutations, the construction can be as follows:
assume that $n$ is even.
If $p>1$, the set $S(n,p)$ has a natural involution $I$ without fixed
points: change the sign of $\pm1$ in permutations.

Then define two subsets of $S(n,p)$ by
\begin{equation}
\begin{split}
S'(n,p)  &= \{\sigma\in S(n,p) | \sigma_{n-1}<-p \text{\ or\ } -1\in\sigma\},\\
S''(n,p) &= \{\sigma\in S(n,p) | \sigma_{n-1}>-p
                                 \text{\ and\ } -1\not\in\sigma\}. \\
\end{split}
\end{equation}
Then $S'(n,p)\cup S''(n,p)$ is $S(n,p)$ and
$S'(n,p) / I(S''(n,p))$ is $S(n,-p)$ up to the sign of the first letter
of each element.
In the special case $p=-1$, both properties still hold, even without the
involution since $S''$ is empty.

Let us illustrate this with the example $n=4$ and $p=2$.
We then have:
\begin{equation}
\begin{split}
S(4,2)=\{&
13\bar42, \bar13\bar42, 14\bar32, \bar14\bar32, \bar31\bar42,
\bar3\bar1\bar42, 3412,
\\ &
34\bar12, \bar3412, \bar34\bar12, \bar41\bar32, \bar4\bar1\bar32, \bar4312,
\bar43\bar12\},
\end{split}
\end{equation}
\begin{equation}
\begin{split}
S'(4,2)=\{&
13\bar42, \bar13\bar42, 14\bar32, \bar14\bar32, \bar31\bar42,
\bar3\bar1\bar42,
\\ &
34\bar12, \bar34\bar12, \bar41\bar32, \bar4\bar1\bar32, \bar43\bar12\},
\end{split}
\end{equation}
\begin{equation}
S''(4,2)=\{3412, \bar3412, \bar4312\},
\end{equation}
so that $S'(4,2)\cup S''(4,2)= S(4,2)$ and
\begin{equation}
\begin{split}
S'(4,2) / I(S''(4,2)) =
\{ &
13\bar42, \bar13\bar42, 14\bar32, \bar14\bar32, \bar31\bar42,
\bar3\bar1\bar42,
\bar41\bar32, \bar4\bar1\bar32\},
\end{split}
\end{equation}
which is $S(4,-2)$ up to the sign of the first letter of each permutation.

\section{Extension to more than two colors}
\label{section7}

\subsection{Augmented alternating permutations}

Let
\begin{equation}
{\bf A} = A^{(0)}\sqcup A^{(1)}\sqcup \dots\sqcup A^{(r-1)}
= A\times C,
\end{equation}
with $C=\{0,\dots,r-1\}$ be an $r$-colored alphabet.
We assume that $A^{(i)}=A\times\{i\}$ is linearly ordered and that
\begin{equation}
A^{(0)} > A^{(1)} > \dots > A^{(r-1)}.
\end{equation}
Colored words can be represented by pairs ${\bf w}=(w,u)$ where
$w\in A^n$ and $u\in C^n$.
We define $r$-colored alternating permutations
${\boldsymbol{\sigma}} = (\sigma,u)$ by the condition
\begin{equation}
  {\boldsymbol\sigma}_1
< {\boldsymbol\sigma}_2
> {\boldsymbol\sigma}_3
< \dots
\end{equation}
hence, as permutations of shape $2^n$ or $2^n1$ as words over $\bf A$.
Let $A_n^{(r)}$ be the set of such permutations.
Their noncommutative generating series in $\FQSym^{(r)}$ is then
\begin{equation}
\sum_{n\geq0} \sum_{{\boldsymbol\sigma}\in A_n^{(r)}}
 \G_{\boldsymbol\sigma}({\bf A})
= {\bf X}({\bf A})
 = \sum_{m\geq0} R_{2^m}({\bf A}) + R_{2^m1}({\bf A}).
\end{equation}
Thus, if we send $A^{(i)}$ to $q_i\E$, the exponential generating function of
the polynomials
\begin{equation}
\alpha_n(q_0,\dots,q_{r-1})
= \sum_{{\boldsymbol\sigma}\in A_n^{(r)}} \prod_{i=1}^n q_{u_i}
\end{equation}
is
\begin{equation}
\alpha(z;q_0,\dots,q_{r-1})
= \frac{\sin(q_0+\dots+q_{r-1})z + 1 }{\cos(q_0+\dots+q_{r-1})z}.
\end{equation}

Iterating the previous constructions, we can define generalized snakes as
colored alternating permutations such that
${\boldsymbol\sigma}_1\in A^{(0)}$,
or, more generally
${\boldsymbol\sigma}_1\in A^{(0)}\sqcup\dots\sqcup A^{(i)}$.
Setting $B_i= A^{(0)}\sqcup\dots\sqcup A^{(i)}$
and $\overline B_i = A^{(i+1)}\sqcup\dots\sqcup A^{(r-1)}$,
we have for these permutations the noncommutative generating series
\begin{equation}
\Y(B_i,\overline B_i) = (\COS + \SIN)(B_i)\ \SEC(B_i|\overline B_i)
\end{equation}
which under the previous specialization yields the exponential generating
series
\begin{equation}
y_i(t;q_0,\dots,q_{r-1})
= \frac{\cos(q_0+\dots+q_i)t + \sin(q_{i+1}+\dots+q_{r-1})t}
       {\cos(q_0+\dots+q_{r-1})t}.
\end{equation}

Setting all the $q_i$ equal to $1$, we recover
sequences~A007286 and~A007289 for $r=3$ and sequence A006873 for $r=4$
of~\cite{Slo}, counting what the authors of~\cite{ER} called
\emph{augmented alternating permutations}.

\subsection{Triangles of alternating permutations with $r$ colors}

One can now count alternating permutations of shapes $2^k$ and $2^k1$ by their
last value. With $r=3$, the following tables present the result:

\medskip
{\scriptsize
\begin{equation}
\label{Ctriang}
\begin{array}{c|ccccccccccccccccc}
n  & {\ol{\ol1}} & {\ol{\ol2}} & {\ol{\ol3}} & {\ol{\ol4}}  & {\ol{\ol5}} &&
     {\ol1} & {\ol2} & {\ol3} & {\ol4} & {\ol5} &&
     1 & 2 & 3 & 4 & 5 \\
\hline
 1 & 1 &   &   &   &   && 1 &   &   &   &   &&
     1 &   &   &   &    \\
 2 & 0 & 1 &   &   &   && 1 & 2 &   &   &   &&
     2 & 3 &   &   &    \\
 3 & 9 & 9 & 8 &   &   && 8 & 7 & 5 &   &   &&
     5 & 3 & 0 &   &    \\
 4 & 0 & 9 &18 &26 &   &&26 &34 &41 &46 &   &&
    46 &51 &54 &54 &    \\
 5 &405&405&396&378&352&&352&326&292&251&205&&
    205&159&108&54 & 0  \\
\end{array}
\end{equation}
}

\begin{proposition}
The table counting alternating permutations with $r$ colors by their last
value is obtained by the following algorithm: first separate the picture by
the column $p=0$ and then compute $r$ triangles. Put $1$ at the top of each
triangle and compute the rest as follows:
fill the second row of all triangles as the sum of the
elements of the first row strictly to their \emph{left}.
Then fill the third row of all triangles as the sum of the
elements of the previous row to their \emph{right}. Compute
all rows successively by reading from left to right and right to left
alternatively.
\end{proposition}

\Proof
Same argument as for Propositions~\ref{tri-alt} and~\ref{tri-sna}.
\qed

Applying the same rules to the construction of three triangles but with only
one $1$ at the top of one triangle gives the following three tables. Note that
this amounts to split the alternating permutations, first by the number of
bars of their \emph{first} element, then, inside the triangle, by their
\emph{last} value.

\medskip
{\scriptsize
\begin{equation}
\label{Ctriangb}
\begin{array}{c|ccccccccccccccccc}
n  & {\ol{\ol1}} & {\ol{\ol2}} & {\ol{\ol3}} & {\ol{\ol4}}  & {\ol{\ol5}} &&
     {\ol1} & {\ol2} & {\ol3} & {\ol4} & {\ol5} &&
     1 & 2 & 3 & 4 & 5 \\
\hline
 1 & 1 &   &   &   &   && 0 &   &   &   &   &&
     0 &   &   &   &    \\
 2 & 0 & 1 &   &   &   && 1 & 1 &   &   &   &&
     1 & 1 &   &   &    \\
 3 & 5 & 5 & 4 &   &   && 4 & 3 & 2 &   &   &&
     2 & 1 & 0 &   &    \\
 4 & 0 & 5 &10 &14 &   &&14 &18 &21 &23 &   &&
    23 &25 &26 &26 &    \\
 5 &205&205&200&190&176&&176&162&144&123&100&&
    100&77 &52 &26 & 0  \\
\end{array}
\end{equation}
}
{\scriptsize
\begin{equation}
\label{Ctriangc}
\begin{array}{c|ccccccccccccccccc}
n  & {\ol{\ol1}} & {\ol{\ol2}} & {\ol{\ol3}} & {\ol{\ol4}}  & {\ol{\ol5}} &&
     {\ol1} & {\ol2} & {\ol3} & {\ol4} & {\ol5} &&
     1 & 2 & 3 & 4 & 5 \\
\hline
 1 & 0 &   &   &   &   && 1 &   &   &   &   &&
     0 &   &   &   &    \\
 2 & 0 & 0 &   &   &   && 0 & 1 &   &   &   &&
     1 & 1 &   &   &    \\
 3 & 3 & 3 & 3 &   &   && 3 & 3 & 2 &   &   &&
     2 & 1 & 0 &   &    \\
 4 & 0 & 3 & 6 & 9 &   && 9 &12 &15 &17 &   &&
    17 &19 &20 &20 &    \\
 5 &147&147&144&138&129&&129&120&108&\ 93\ &\ 76\ &&
   \,\,\, 76\,\,\ &59 &40 &20 & 0  \\
\end{array}
\end{equation}
}
{\scriptsize
\begin{equation}
\label{Ctriangd}
\begin{array}{c|ccccccccccccccccc}
n  & {\ol{\ol1}} & {\ol{\ol2}} & {\ol{\ol3}} & {\ol{\ol4}}  & {\ol{\ol5}} &&
     {\ol1} & {\ol2} & {\ol3} & {\ol4} & {\ol5} &&
     1 & 2 & 3 & 4 & 5 \\
\hline
 1 & 0 &   &   &   &   && 0 &   &   &   &   &&
     1 &   &   &   &    \\
 2 & 0 & 0 &   &   &   && 0 & 0 &   &   &   &&
     0 & 1 &   &   &    \\
 3 & 1 & 1 & 1 &   &   && 1 & 1 & 1 &   &   &&
     1 & 1 & 0 &   &    \\
 4 & 0 & 1 & 2 & 3 &   && 3 & 4 & 5 & 6 &   &&
     6 & 7 & 8 & 8 &    \\
 5 &53 &53 &52 &50 &47 &&47 &44 &40 &35 &29 &&
    29 &23 &16 & 8 & 0  \\
\end{array}
\end{equation}
}

Arnol'd \cite{Arnold} has found remarkable arithmetical properties
of the Euler-Bernoulli triangles.
The study of the properties of these new triangles remains to be done.

\section{Snakes of type $D$}
\label{sectionD}

\subsection{The triangle of type $D$ snakes}

The Springer numbers of type $D$ (Sequence A007836 of~\cite{Slo}) are given by
exactly the same  process as for type $B$ Springer numbers, but starting with
$0$ at the top of the left triangle and $1$ at the top of the right triangle.
Since all operations computing the rows of the triangles are linear in the
first entries, we have in particular that the sum  of the number of snakes of
type $D_n$ and  the number of snakes on type $B_n$ is equal to the number of
alternating permutations of type $B$.

We have even more information related to the triangles: both $B$ and
$D$ triangles can be computed by taking the difference between the triangle of
Equation~(\ref{AltBtriang}) and of Equation~(\ref{EBtriang}).
We  obtain

\begin{equation}
\label{Dtriang}
\begin{array}{c|ccccccccccccc}
n/p&-6 &-5 &-4 &-3 &-2 &-1 &   & 1 & 2 & 3 & 4 &5  &6 \\
\hline
 1 &   &   &   &   &   & 0 &   & 1 &   &   &   &   &  \\
 2 &   &   &   &   & 0 & 0 &   & 0 & 1 &   &   &   &  \\
 3 &   &   &   & 1 & 1 & 1 &   & 1 & 1 & 0 &   &   &  \\
 4 &   &   & 0 & 1 & 2 & 3 &   & 3 & 4 & 5 & 5 &   &  \\
 5 &   &23 &23 &22 &20 &17 &   &17 &14 &10 & 5 & 0 &  \\
 6 & 0 &23 &46 &68 &88 &105&   &105&122&136&146&151&151\\
\end{array}
\end{equation}

\medskip
\subsubsection{Snakes of type $D$}

From our other sets having the same cardinality as type $B$ snakes, we can
deduce combinatorial objects having same cardinality as type $D$ snakes by
taking the complement in the alternating permutations of type $B$.

Since the generating series of type $B$ alternating permutations is
\begin{equation}
\X =\sum_{m\ge 0}(R_{(2^m)}^\# + R_{(2^m1)}^\#).
\end{equation}
and the generating series of type $B$ snakes defined in
Section~\ref{snakesssens} is
\begin{equation}
\Y = 1+\sum_{n\ge 0}(\tilde R_{02^n1} +\tilde R_{02^{n+1}}),
\end{equation}
we easily get one definition of the generating series type $D$ snakes:
\begin{equation}
\label{dsnakes}
\D = \X - \Y = \sum_{n\ge 0}(\tilde R_{2^n1} +\tilde R_{2^{n+1}}).
\end{equation}

In other words, our first sort of type $D$ snakes corresponds to
permutations of ribbon shape $2^n1$ or $2^n$ whose first letter is positive.

\medskip
{\scriptsize
Here are these elements for $n\leq4$:
\begin{equation}
\quad
12, 
\quad
12\bar3 , 132 , 231 , 13\bar2 , 23\bar1
\end{equation}
\begin{equation}
\begin{split}
&
\bar12\bar34, 12\bar43, 12\bar4\bar3, 1324, 13\bar24,
13\bar42, 13\bar4\bar2, 1423, 14\bar23, 14\bar32,
14\bar3\bar2, 2314,
\\ &
23\bar14, 23\bar41, 23\bar4\bar1,
2413, 24\bar13, 24\bar31, 24\bar3\bar1, 3412,
34\bar12, 34\bar21, 34\bar2\bar1.
\end{split}
\end{equation}
}

Since both alternating permutations and snakes of type $B$ can be interpreted
as solutions of a differential equation and a fixed point solution involving
the same bilinear form, one then concludes that these snakes of type $D$
satisfy
\begin{equation}
d \D = 1 + \D \X,
\end{equation}
and
\begin{equation}
\label{fixD}
\D = \G_1 + \B(\D,\X).
\end{equation}
The iteration of~(\ref{fixD})  brings up a solution close to ~\eqref{X2cbt}
and~\eqref{Y2cbt}:
\begin{equation}
\D = \sum_{T\in \CBT} B_T(\G_0,\G_1),
\end{equation}
where, for a tree $T$, $B_T(a,b)$ is the result of the evaluation of all
expressions formed by labeling by $a$ or $b$ the leaves of $T$ and by $B$ its
internal nodes. Note that in this case, the first leaf needs to have label 
$b$.

\subsubsection{The usual snakes of type $D$}

The previous type $D$ snakes are not satisfactory since, even if they fit into
the desired triangle, they do not belong to $D_n$.
The classical snakes of type $D$ of Arnol'd belong to $D_n$ and are easily
defined: select among permutations of ribbon shape $12^n$ and $12^n1$ the
elements with an even number of negative signs and such that
$\sigma_1+\sigma_2<0$.

\medskip
{\scriptsize
One then gets the following elements for $n\leq4$:
\begin{equation}
1\quad
\bar1\bar2, 
\quad
\bar1\bar23, 1\bar3\bar2, \bar1\bar32, 2\bar3\bar 1, \bar2\bar31,
\end{equation}
\begin{equation}
\begin{split}
&
1\bar23\bar4, 1\bar24\bar3, \bar1\bar243, 1\bar32\bar4,
\bar1\bar3\bar2\bar4, 1\bar34\bar2,
\bar1\bar342, 1\bar42\bar3, \bar1\bar4\bar2\bar3, 1\bar43\bar2,
\bar1\bar432, 2\bar31\bar4,
\\ &
\bar2\bar3\bar1\bar4, 2\bar34\bar1, \bar2\bar341, 2\bar41\bar3,
\bar2\bar4\bar1\bar3, 2\bar43\bar1,
\bar2\bar431, 3\bar41\bar2, \bar3\bar4\bar1\bar2, 3\bar42\bar1,
\bar3\bar421.
\end{split}
\end{equation}
}

\medskip
It is easy to go from these last elements to the other type $D$ snakes:
change all values into their opposite and then change the first element $s$ to
$-|s|$. Conversely, change (resp. do not change) the sign of the first element
depending whether it is not (resp. it is) in $D_n$ and then change all values
into their opposite.

\section{Tables}
\label{ann-tilde}

Here follow the tables of the maps $\ol A \mapsto tA$ from $\BSym$ to $\NCSF$.

All tables represent in columns the image of ribbons indexed by type-$B$
compositions, where the first half begins with a $0$
and the other half does not.
So, with $N=3$, compositions are in the following order:
\begin{equation}
[0,3],\ [0,2,1],\ [0,1,2],\ [0,1,1,1],\
[  3],\ [  2,1],\ [  1,2],\ [  1,1,1].
\end{equation}
Note that the zero entries have been represented by dots to enhance
readability.


\begin{figure}[ht]
\begin{equation}
\left(
\begin{array}{cccc}
-t^2 & t^2 & 1 & -1 \\
t^2 + t & . & . & t + 1
\end{array}
\right)
\quad
\left(
\begin{array}{cccccccc}
t^3 & -t^3      & -t^3        & t^3 & 1     & -1 & -1 & 1 \\
- t^3\! -\! t^2 & .           & t^3\!+\!t^2 & . & . & t\! +\! 1
           & . & - t\! -\! 1 \\
t\! -\! t^3     & t^3 \!-\! t & . & . & . & . & 1\! -\! t^2 & t^2\! -\! 1 \\
t^3\! +\! t^2   & t^2\! +\! t & . & . & . & . & t^2\! +\! t & t\! +\! 1 
\end{array}
\right)
\end{equation}
\caption{
\label{fig-S-23}
Matrices of $\tilde R_I(A,tA)$ on the $S$ basis for $n=2$, $3$.}
\end{figure}


\begin{figure}[ht]
\begin{equation}
\left(
\begin{array}{cccc}
t^2 & -t^2 & -1 & 1 \\
t & t^2 & 1 & t
\end{array}
\right)
\quad
\left(
\begin{array}{cccccccc}
t^3 & -t^3 & -t^3 & t^3 & 1 & -1 & -1 & 1 \\
t^2 & t^3 & -t^2 & -t^3 & -1 & -t & 1 & t \\
-t & t & t^3 & -t^3 & -1 & 1 & t^2 & -t^2 \\
t & t^2 & t^2 & t^3 & 1 & t & t & t^2
\end{array}
\right)
\end{equation}
\caption{
\label{fig-L-23}
Matrices of $\tilde R_I(A,tA)$ on the $\Lambda$ basis for $n=2$, $3$.}
\end{figure}


\begin{figure}[ht]
\begin{equation}
\left(
\begin{array}{cccc}
t & t^2 & 1 & t \\
t^2 + t & . & . & t + 1
\end{array}
\right)
\quad
\left(
\begin{array}{cccccccc}
t & t^2 & t^2 & t^3 & 1 & t & t & t^2 \\
. & t^2 + t & t^3 + t^2 & . & . & t + 1 & t^2 + t & . \\
t^2 + t & t^3 + t^2 & . & . & . & . & t + 1 & t^2 + t \\
t^3 + t^2 & t^2 + t & . & . & . & . & t^2 + t & t + 1 
\end{array}
\right)
\end{equation}
\caption{
\label{fig-R-23}
Matrices of $\tilde R_I(A,tA)$ on the $R$ basis for $n=2$, $3$.}
\end{figure}


\end{document}